\documentstyle[12pt]{article}

\newcommand{\mn}{\sl}
\def\afterthmseparator{}
\makeatletter
\renewcommand{\@begintheorem}[2]{\trivlist
      \item[\hskip \labelsep{\bf #1\ #2\unskip\afterthmseparator}]\mn}
\renewcommand{\@opargbegintheorem}[3]{\trivlist
      \item[\hskip \labelsep{\bf #1\ #2\ (#3)\unskip\afterthmseparator}]\mn}
\makeatother
\newtheorem{theorem}{Theorem}[section]
\newtheorem{lemma}[theorem]{Lemma}
\newtheorem{corollary}[theorem]{Corollary}
\newtheorem{proposition}[theorem]{Proposition}
\newtheorem{rem}[theorem]{Remark}
\newenvironment{remark}{\renewcommand{\mn}{\rm} \begin{rem}}{\end{rem}}
\newtheorem{probl}[theorem]{Problem}

\newtheorem{df}[theorem]{Definition}
\newenvironment{definition}{\renewcommand{\mn}{\rm} \begin{df}}{\end{df}}
\newtheorem{exmpl}[theorem]{Example}
\newenvironment{example}{\renewcommand{\mn}{\rm} \begin{exmpl}}{\end{exmpl}}
\newcounter{oq}
\newcommand{\que}{\refstepcounter{oq}\par{\bf \theoq.}~}

\newcommand{\bull}{\mbox{$\;\;\;$\vrule height .9ex width .8ex depth -.1ex}}
\newcommand{\qed}{\mbox{$\;\;\;\Box$}}
\newcommand{\noproof}{\unskip\nobreak\hfill \bull}
\newenvironment{proof}{\par\smallbreak\noindent{\sl Proof.~}}%
{\unskip\nobreak\hfill \bull \par\medbreak}
\newenvironment{proofof}[1]{\par\smallbreak\noindent{\sl Proof~of~#1.~}}%
{\unskip\nobreak\hfill \bull \par\medbreak}
\newenvironment{subproof}{\par\noindent{\it Proof of Claim.~}}%
{\qed \par\smallbreak}
\newcounter{claim}[theorem]
\renewcommand{\theclaim}{\thetheorem.\arabic{claim}}
\newenvironment{claim}{\refstepcounter{claim}%
\par\medskip\par\noindent{\it Claim~\theclaim.~}~\rm}%
{\par\smallskip\par}
\newcommand{\case}[2]{\par{\it Case #1:}\/ #2.}

\newcommand{\hide}[1]{}
\newcommand{\refeq}[1]{(\ref{#1})}
\newcommand{\function}[2]{:#1 \rightarrow #2}
\newcommand{\of}[1]{\left( #1 \right)}
\newcommand{\setdef}[2]{\left\{\hspace{0.5mm}#1:\hspace{0.5mm} #2\right\}}

\newcommand{\ga}{\alpha}

\newcommand{\gs}{\sigma}
\newcommand{\Gg}{\Gamma}
\newcommand{\And}{\wedge}
\newcommand{\Or}{\vee}

\newcommand{\EF}{Ehrenfeucht}
\newcommand{\WL}{Weisfeiler-Lehman}
\newcommand{\Wl}{\mbox{\rm WL}}
\newcommand{\game}{\mbox{\sc Ehr}}
\newcommand{\Ii}[1]{\mbox{\rm\=D} (#1)}
\newcommand{\I}[1]{\Ii{#1}}
\newcommand{\II}[2]{\mbox{\rm\=D}_{#1} (#2)}
\newcommand{\Vv}[1]{\mbox{\=L} (#1)}
\newcommand{\Cc}[1]{\mbox{\=C} (#1)}
\newcommand{\D}[1]{\mbox{\rm D} (#1)}
\newcommand{\DD}[2]{\mbox{\rm D}_{#1} (#2)}
\newcommand{\V}[1]{\mbox{\rm L} (#1)}
\newcommand{\C}[1]{\mbox{\rm C} (#1)}
\newcommand{\cl}{{\cal C}}
\newcommand{\cls}{{\cal S}}
\newcommand{\Nest}{\mathop{\rm Nest}\nolimits}
\newcommand{\nest}[1]{\Nest(#1)}
\newcommand{\Qr}{\mathop{\rm qr}\nolimits}
\newcommand{\qr}[1]{\Qr(#1)}
\newcommand{\Conn}{\mathop{\rm conn}\nolimits}
\newcommand{\conn}[1]{\Conn(#1)}
\newcommand{\Alt}{\mathop{\rm alt}\nolimits}
\newcommand{\alt}[1]{\Alt(#1)}
\newcommand{\compl}[1]{\overline{#1}}
\newcommand{\add}{\oplus}

\newcommand{\xeq}{\mathbin{{\equiv}_X}}

\newcommand{\xyxeq}{\mathbin{{\equiv}_{X\cup Y(X)}}}
\newcommand{\xxeq}{\mathbin{{\equiv}_{X'}}}
\newcommand{\phieq}{\mathbin{{\equiv}_\phi}}
\newcommand{\phistareq}{\mathbin{{\equiv}_{\phi^*}}}
\newcommand{\notphistar}{\mathbin{{\not\equiv}_{\phi^*}}}
\newcommand{\notxyxk}{\mathbin{{\not\equiv}_{X_{k-1}\cup Y(X_{k-1})}}}
\newcommand{\calC}{{\cal C}}
\newcommand{\calD}{{\cal D}}
\newcommand{\dist}{\mathop{\rm dist}\nolimits}
\newcommand{\baru}{{\bar u}}
\newcommand{\barv}{{\bar v}}
\newcommand{\isotype}[1]{[#1]}
\newcommand{\wl}[4]{W_{#1}^{#2,#3}(#4)}
\newcommand{\wlcol}[3]{W_{#1}^{#2,#3}}
\newcommand{\vocab}{{\cal L}}

\newif\ifnotesw\noteswtrue

\title{
The First Order Definability of Graphs:\\
Upper Bounds for Quantifier Rank}

\author{Oleg Pikhurko\thanks{%
Department of Mathematical Sciences,
Carnegie Mellon University, Pittsburgh, PA 15213-3890
Web: {\tt http://www.math.cmu.edu/\~{}pikhurko/}}
\quad
Helmut Veith\thanks{%
Institut f\"ur Informationssysteme, Technische Universit\"at Wien,
Favoritenstr.~9, A-1040 Wien, Austria.
Supported by the European Community Research Training Network
``Games and Automata for Synthesis and Validation'' (GAMES)
and by the Austrian Science Fund Project Z29-INF.
E-mail: {\tt veith@dbai.tuwien.ac.at}}
\quad
Oleg Verbitsky\thanks{%
Department of Mechanics \& Mathematics, Kyiv University, Ukraine.
Research was done in part while visiting the
Institut f\"ur Informationssysteme at the Technische Universit\"at Wien,
supported from Austrian Science Foundation grant Z29-INF.
E-mail: {\tt oleg@ov.litech.net}}}

\date{}

\begin{document}
\maketitle

\begin{abstract}
We say that a first order formula $\Phi$ distinguishes a graph $G$ from
another graph $G'$ if $\Phi$ is true on $G$ and false on $G'$.
Provided $G$ and $G'$ are non-isomorphic,
let $\D{G,G'}$ denote the minimal quantifier rank of a such formula.
Let $n$ denote the order of $G$.
We prove that, if $G'$ has the same order, then
$\D{G,G'}\le(n+3)/2$. This bound is tight up to an additive
constant of 1.
Furthermore, we prove that non-isomorphic $G$ and $G'$ of order $n$
are distinguishable by an existential formula of quantifier rank
at most $(n+5)/2$.
As a consequence of the first result, we obtain an upper bound of $(n+1)/2$
for the optimum dimension of the \WL\/ graph canonization
algorithm, whose worst case value is known to be linear in~$n$.

We say that a first order formula $\Phi$ defines a graph $G$ if $\Phi$
distinguishes $G$ from every non-isomorphic graph $G'$.
Let $\D G$ be the minimal quantifier rank of a formula defining $G$.
As it is well known, $\D G\le n+1$ and this bound
is generally best possible. Nevertheless, we here show that there is
a class $\cl$ of graphs of simple, easily recognizable structure such that
\begin{itemize}
\item
$\D G\le (n+5)/2$ with the exception of all graphs in~$\cl$;
\item
if $G\in\cl$, then it is easy to compute the exact value of~$\D G$.
\end{itemize}
Moreover, the defining formulas in this result have only one quantifier
alternation.
The bound for $\D G$ can be improved for graphs with bounded vertex degrees:
For each $d\ge2$ there is a constant $c_d<1/2$ such that
$\D G\le c_d n+O(d^2)$ for any graph $G$ with no isolated vertices
and edges whose maximum degree is~$d$.

Finally, we extend our results over directed graphs, more
generally, over arbitrary structures with maximum relation arity 2,
and over $k$-uniform hypergraphs.
\end{abstract}

\tableofcontents

\clearpage

\section{Introduction}\label{s:intro}

{}From the logical point of view, a graph $G$ is a structure with a single
anti-reflexive and symmetric binary predicate $E$ for the adjacency
relation of $G$. Every closed first order formula $\Phi$ with predicate
symbols $E$ and $=$ is either true or false on $G$.
Given two non-isomorphic graphs $G$ and $G'$, we say that $\Phi$
{\em distinguishes $G$ from $G'$} if $\Phi$ is true on $G$ but false on
$G'$. The {\em quantifier rank\/} of $\Phi$ is the maximum
number of nested quantifiers in this formula (see Section \ref{s:prel} for
formal definitions). Let $\D{G,G'}$ denote the minimum quantifier rank of
a formula distinguishing $G$ from $G'$. The number $\D{G,G'}$
is symmetric with respect to $G$ and $G'$ because,
if $\Phi$ distinguishes $G$ from $G'$, then $\neg\Phi$
distinguishes $G'$ from $G$ and has the same quantifier rank.

The {\em order\/} of a graph is the number of its vertices.
Throughout the paper the letter $n$ will denote the order of
a graph $G$.

The first question we address is how large $\D{G,G'}$ can be
over non-isomorphic graphs $G$ and $G'$ of the same order $n$.
For any $n$ it is not hard to find a such pair with
$$
\D{G,G'}\ge(n+1)/2
$$
(see Example \ref{ex:lower}).
On the other hand, there is an obvious general upper bound
$$
\D{G,G'}\le n.
$$
Indeed, a graph $G$ with vertex set $V(G)=\{1,\ldots,n\}$
and edge set $E(G)$ is distinguished from any non-isomorphic $G'$
of order $n$ by the formula
\begin{equation}\label{eq:dist}
\exists x_1\ldots\exists x_n\of{
\bigwedge_{i\ne j} \neg (x_i=x_j) \And
\bigwedge_{\{i,j\}\in E(G)} E(x_i,x_j) \And
\bigwedge_{\{i,j\}\notin E(G)} \neg E(x_i,x_j)}.
\end{equation}
It seems that no better upper bound has been reported in the literature
so far. Here we prove a nearly best possible bound.

\begin{theorem}\label{thm:d}
If $G$ and $G'$ are non-isomorphic graphs both of order $n$, then
$$
\D{G,G'}\le(n+3)/2.
$$
\end{theorem}

It is worth noting that the distinguishing formulas resulting from
our proof of Theorem \ref{thm:d}
have a rather restricted logical structure. We say that a first order formula
$\Phi$ is in the {\em negation normal form\/} if the connective $\neg$
occurs in $\Phi$ only in front of atomic subformulas. If $\Phi$ is
such a formula, its {\em alternation number\/} is the maximum number
of alternations of $\exists$ and $\forall$ in a sequence of nested
quantifiers of $\Phi$.
The proof of Theorem \ref{thm:d} produces distinguishing formulas
in the negation normal form whose alternation number is at most 1.
We are able to prove Theorem \ref{thm:d} even with alternation number 0
but with a little weaker bound $\D{G,G'}\le(n+5)/2$. The proof actually
produces either an existential or a universal distinguishing formula.
An {\em existential\/} or {\em universal\/}
formula is a formula in the negation normal form whose all quantifiers
are of only one sort, existential or universal respectively.

Our proof of these results is based on the well-known combinatorial
characterization of $\D{G,G'}$ in terms of the \EF\/
game on $G$ and $G'$ \cite{Ehr} (Fra\"\i ss\'e \cite{Fra}
suggested an essentially equivalent characterization in terms of
partial isomorphisms between $G$ and $G'$).
In this setting, $\D{G,G'}$ is equal to the length of the game
under the condition that the players play optimally.
Thus Theorem \ref{thm:d} says that the \EF\/ game on non-isomorphic
graphs of the same order $n$ can be won within at most $(n+3)/2$ rounds.

Theorem \ref{thm:d} has consequences for the complexity analysis
of the \WL\/ algorithm for graph isomorphism testing. This algorithm
has been studied since the seventies (see e.g.\ \cite{Bab,CFI}).
An important combinatorial parameter of the algorithm,
occurring in the known bounds on the running time, is its {\em dimension}.
Denote the optimum dimension for input graphs $G$ and $G'$ by
$\Wl(G,G')$ (see Section \ref{s:wl} for a detailed exposition).

Cai, F\"urer, and Immerman \cite{CFI} come up with a remarkable
construction of non-isomorphic $G$ and $G'$ of the same order
for which $\Wl(G,G')=\Omega(n)$. Though they do not specify
the constant hidden in the $\Omega$-notation, a simple analysis of
their proof, that we make in Section \ref{ss:cfi},
gives at least
\begin{equation}\label{eq:cfi}
\Wl(G,G') > 0.00465\,n.
\end{equation}
In the same paper, the authors give a logical characterization of
$\Wl(G,G')$ which readily implies the relation
\begin{equation}\label{eq:cfichar}
\Wl(G,G')\le\D{G,G'}-1.
\end{equation}
Thus, our Theorem \ref{thm:d} establishes an upper bound
\begin{equation}\label{eq:wlupper}
\Wl(G,G')\le 0.5\,n+0.5.
\end{equation}
We have to remark that this bound for the dimension does not imply
any good bound for the worst case running time. In fact, the lower bound
\refeq{eq:cfi} shows that the \WL\/ algorithm hardly can be practical
in the worst case. However, we believe that \refeq{eq:wlupper} shows
an interesting combinatorial property of the algorithm previously
never observed.

Providing upper bounds for the length of the \EF\/ game and for
the dimension of the \WL\/ algorithm, Theorem \ref{thm:d} is
therefore rather meaningful from combinatorial and computer
science point of view. Let us now focus on logical motivations
in the scope of finite model theory.

We say that a formula $\Phi$ {\em identifies\/}
a graph $G$ (up to an isomorphism in the class of graphs of the same order)
if $\Phi$ distinguishes $G$ from
any other non-isomorphic graph $G'$ of the same order.
Let $\I G$ denote the minimum quantifier rank of a such formula.
Note that, if $G$
is distinguished from $G'$ by formula $\Phi_{G'}$, then $G$ is identified
by the conjunction $\bigwedge_{G'}\Phi_{G'}$ over all non-isomorphic $G'$
of the same order. It easily follows that $\I G=\max \D{G,G'}$ over
all such $G'$. In these terms
Theorem \ref{thm:d} reads $\I G\le(n+3)/2$.

However, from the point of view of finite model theory it is more
natural to address defining rather than distinguishing formulas.
A first order formula $\Phi$ {\em defines\/} a graph $G$ if $\Phi$
distinguishes $G$ from all non-isomorphic graphs, regardless of their
order. Let $\D G$ be the minimum quantifier rank of a
formula defining $G$. As it is well known,
$$
\D G\le n+1
$$
for every $G$ (we have to append $\forall x_{n+1}$ to the quantifier
prefix of \refeq{eq:dist} and say that any vertex $x_{n+1}$ is
equal to one of $x_1,\ldots,x_n$). This bound cannot be generally
improved because, for example, no formula of quantifier rank $n$
can distinguish between two complete graphs of orders $n$ and $n+1$.
Nevertheless, our next aim is to suggest a bound for $\D G$ similar to
Theorem \ref{thm:d} and explicitly describe all the exceptions.
We are able to prove a dichotomy result: Either $\D G\le n/2+O(1)$
or else $G$ has a simple structure and, moreover, in the latter case
$\D G$ is efficiently computable.

\begin{theorem}\label{thm:def}
There is an efficiently recognizable class of graphs $\cl$ such that
\begin{itemize}
\item
$\D G\le (n+5)/2$ with the exception of all graphs in~$\cl$;
\item
if $G\in\cl$, then the exact value of $\D G$ is efficiently computable.
\end{itemize}
Moreover, every graph $G$ admits a defining formula whose alternation number
is 1 and whose quantifier rank is as small as possible if $G\in\cl$ and
does not exceed $(n+5)/2$ if $G\notin\cl$.
\end{theorem}

Referring to the efficiency here, we mean the time $O(n^2\log n)$ on a
random access machine with input graphs given by their adjacency matrices.
T.~\L uczak \cite{Luc} poses a question if $\D G$ is computable.
It seems plausible that it is not. The simulation of undecidable
problems by finite structures, as in the Trakhtenbrot theorem \cite{Tra}
(cf.\ \cite{Vau} and \cite[theorem 8.2.1]{Spe}),
can be considered an evidence in favor of this
hypothesis. In this respect, Theorem \ref{thm:def} provides us with
a non-trivial computable upper bound for $\D G$.

The bound of Theorem \ref{thm:def} can be improved for graphs
with bounded vertex degree:
For each $d\ge2$ there is a constant $c_d<1/2$ such that
$\D G\le c_d n+O(d^2)$ for any graph $G$ with no isolated vertices
and edges whose maximum degree is $d$.
We do not try to find the best possible $c_d$ being content with
a constant strictly less than $1/2$. Note that no sublinear bound
is possible here. It is easy to show, for example, that
$\D{G_n}\ge n/(d+1)$ for $n=m(d+1)$ and $G_n$ being the vertex disjoint
union of $m$ copies of the complete graph on $d+1$ vertices.
Moreover, no sublinear upper bound is possible even for connected graphs,
as follows from the Cai-F\"urer-Immerman bounds \refeq{eq:cfi} and
\refeq{eq:cfichar} and because $G$ and $G'$ in their construction
are connected graphs of bounded degree (see Section~\ref{ss:cfi}).

With minor efforts, Theorems \ref{thm:d} and \ref{thm:def}
carry over to directed graphs and, more generally, to arbitrary
relational structures with maximum arity 2. In combinatorial terms,
the latter are directed graphs endowed with colorings of the vertex
set and the edge set.

Finally, we prove somewhat weaker analogs of
Theorems \ref{thm:d} and \ref{thm:def} for $k$-uniform hypergraphs.
The upper bound we obtain here is $(1-1/k)n+2k-1$. It remains open
if this bound is tight for $k\ge3$ since the only lower bound we
know for any $k$ is $(n+1)/2$.

\subsubsection*{Previous work}
The graph identification is studied in \cite{IKo,ILa,CFI,Gro1,Gro2,GMa}
in aspects relevant to computer science. The main focus of this line
of research is on the minimum number of variables used in an identifying
formula, where formulas are in the first order language enriched by
counting quantifiers. A {\em counting quantifier\/}
$\exists^{\ge m}x\Psi$ means that there are at least $m$
vertices $x$ for which the statement $\Psi$ holds.
Let $\V{G,G'}$ denote the minimum number of variables in a
formula distinguishing non-isomorphic graphs $G$ and $G'$
(different occurrences of the same variable are not counted).
Let $\C{G,G'}$ be the analog of this number
for the logic with counting quantifiers.
Since every formula of quantifier rank $r$ can be rewritten
in equivalent form using only $r$ variables, we have
$$
\C{G,G'}\le\V{G,G'}\le\D{G,G'}.
$$
Asuming that $G$ and $G'$ are of the same order, our Theorem \ref{thm:d}
establishes an upper bound of $(n+3)/2$ for the whole hierarchy.
The aforementioned lower bound by Cai, F\"urer, and Immerman \cite{CFI}
is actually $\C{G,G'}=\Omega(n)$, for infinitely many non-isomorphic
$G$ and $G'$. Their characterization of the optimum dimension
of the \WL\/ algorithm is $\Wl(G,G')=\C{G,G'}-1$.

Let $\Vv G$ (resp.\ $\Cc G$) be the minimum number of variables
in a formula (resp.\ with counting quantifiers) identifying a graph $G$.
By the aforementioned relationship between distinguishing and
identifying formulas, we have $\Vv G=\max \V{G,G'}$ over
all $G'$ of order $n$ non-isomorphic with $G$. Similarly,
$\Cc G=\max \C{G,G'}$ over $G'$ non-isomorphic with $G$,
where $G'$ may be of any order since graphs of different orders
are easily distinguishable in the logic with counting quantifiers.
Thus, $\Cc G$ is equal to the minimum number of variables in a formula
with counting quantifiers that {\em defines\/} $G$.

Define $\I n$, $\Vv n$, and $\Cc n$ to be the maximum possible values
of $\I G$, $\Vv G$, and $\Cc G$, respectively, over graphs of order $n$.
We now know that
$$
0.00465\,n<\Cc n\le\Vv n\le\I n\le 0.5\,n+1.5,
$$
where the upper bound is given by Theorem \ref{thm:d} and the
lower bound is actually the Cai-F\"urer-Immerman bound \refeq{eq:cfi}.
The values of $\Vv n$ and $\I n$ are at most 1 apart from the upper
bound. An interesting open question is where in this range $\Cc n$
is located.

It is known that $\Cc G=2$ for almost all graphs of a given order,
$\Cc G=3$ for almost all regular graphs of a given order and degree,
$\Cc G=2$ for all trees \cite{BES,BKu,Kuc,ILa},
$\Cc G=O(g)$ for graphs of genus $g$ \cite{Gro1,Gro2}, and
$\Cc G\le k+2$ for graphs of tree-width $k$ \cite{GMa}.
For strongly regular graphs it holds $\Cc G=O(\sqrt n\log n)$,
where $n$ denotes the order of $G$ \cite{Bab2}.
If $G$ has separator of size $O(n^\delta)$, $0<\delta<1$,
then $\Cc G=O(n^\delta)$ \cite{CFI}. This result applies, in particular,
to classes of graphs with excluded minors, that have separators of size
$O(\sqrt n)$~\cite{AST}.

Estimation of $\D{G,G'}$ is an interesting research problem not only for
graphs but also for any class of structures. The case of words is
considered in \cite{SJo} in the context of Zero-One Laws of first order
logic, where $\D{W,W'}$ is estimated for $W$ and $W'$ being independent
random binary words of length~$n$.

\subsubsection*{Subsequent work}

In \cite{Ver} it is shown that $\D G=O(\log n)$ if $G$ is a tree
of bounded degree or a Hamiltonian outerplanar graph. This upper bound
complements the popular lower bounds
$\D{P_n,P_{n+1}}>\log_2n-3$
(e.g.\ \cite[Theorem 2.1.3]{Spe}) and
$\D{C_n,C_{n+1}}>\log_2n$
(e.g.\ \cite[Example 2.3.8]{EFl}), where
$P_n$ is the path and $C_n$ is the cycle on $n$ vertices.

As already mentioned, in the present paper we manage to extend
Theorem \ref{thm:d} from graphs to arbitrary structures with
maximum relation arity 2 and to $k$-uniform hypergraphs.
Extension to arbitrary structures $G$ and $G'$ with maximum relation
arity $k$ seems a much more subtle problem. It is done in \cite{PVe}
with bound $(1-\frac1{2k})n+k^2-k+2$, even if we restrict the
alternation number of a distinguishing formula to 1. The same bound
is established for $\D G$ provided no transposition of two elements
of $G$ is an automorphism of the structure.

Theorem \ref{thm:d} gives a worst case upper bound on $\D{G,G'}$
for graphs of the same order. In \cite{KPSV} the average case is
analyzed and it is proved that, for random independent $G$ and $G'$
on $n$ vertices, with probability $1-o(1)$ we have
$\D{G,G'}=\log_2n(1+o(1))$.
Moreover, $\D G$ for a random graph $G$ is determined
with high precision: With probability $1-o(1)$,
$$
\log_2n-2\log_2\log_2n\le \D G
\le\log_2n-\log_2\log_2n+\log_2\log_2\log_2n+O(1).
$$
The upper bound here holds even if the alternation number of
defining formulas is restricted to 1.
A logarithmic upper bound is also proved with the alternation number
restricted to 0. Together with Theorem \ref{thm:D0}
in the present paper, this shows that bounding the alternation number
does not affect much the maximal and the average values of $\D{G,G'}$
for graphs of the same order.

In \cite{PSV} the ``best case'' behavior of $\D G$ is investigated.
Namely, let $g(n)=\min\D G$ over graphs of order $n$.
It is not hard to see that $g(n)\to\infty$ as $n\to\infty$ but it is
not so clear how fast or slowly $g(n)$ grows. In \cite{PSV}
it is proved that
$g(n)$ can be so small if compared to $n$ that the gap between the
two numbers cannot be estimated by any computable function, i.e.,
there is no recursive function $f$ such that $n\le f(g(n))$ for all $n$.
However, if we ``regularize'' $g(n)$ by considering
$\bar g(n)=\max_{m\le n}g(m)$, we have
$\bar g(n)=(1+o(1))\log^* n$, where $\log^*n$ is the smallest number of
iterations of the logarithm base 2 that suffices to decrease $n$
below 1. This result is proved even under the restriction of
the alternation number of a defining formula to a constant.
Under the strongest restriction of the alternation number to 0,
an infinite family of graphs with $\D G\le2\log^*n+O(1)$
is constructed.

\subsubsection*{Organization of the paper}

The paper is organized as follows. Section \ref{s:prel} contains
the relevant definitions from graph theory and logic as well as the
basic facts on the \EF\/ games. In Sections \ref{s:dist} and \ref{s:def}
we prove Theorems \ref{thm:d} and \ref{thm:def}, restated there as Theorems
\ref{thm:d1} and \ref{thm:def1}. In Section \ref{s:d0} we prove
a variant of Theorem \ref{thm:d} for distinguishing formulas whose
alternation number is restricted to 0. The case of graphs with
bounded vertex degree is considered in Section \ref{s:bdeg}.
Section \ref{s:wl} gives an exposition of the \WL\/ algorithm
and applies Theorem \ref{thm:d} for its analysis. This section is
mostly expository and, in particular, includes computation of
an explicit constant in the Cai-F\"urer-Immerman bound \refeq{eq:cfi}.
Section \ref{s:digraphs} discusses the extensions
of Theorems \ref{thm:d} and \ref{thm:def} to directed graphs and
relational structures of maximum arity 2.
Section \ref{s:hyper} is devoted to $k$-uniform hypergraphs.
Section \ref{s:open} lists open problems.

\section{Preliminaries}\label{s:prel}

\subsection{Structures}

\subsubsection{Graphs}

Given a graph $G$, we denote its vertex set by $V(G)$ and its edge set
by $E(G)$. The order of $G$ will be sometimes denoted by $|G|$, that is,
$|G|=|V(G)|$. The {\em neighborhood\/} of a vertex $v$ consists of
all vertices adjacent to $v$ and is denoted by $\Gg(v)$.

The {\em complement\/} of $G$, denoted by $\compl G$, is
the graph on the same vertex set $V(G)$ with all those edges that are
not in $E(G)$.
Given $G$ and $G'$ with disjoint vertex sets, we define
the {\em sum} (or {\em disjoint union}) $G\sqcup G'$ to be the graph
with vertex set $V(G)\cup V(G')$ and edge set $E(G)\cup E(G')$.

A set $S\subseteq V(G)$ is called {\em independent\/} (or {\em stable\/})
if it contains no pair of adjacent vertices. $S$ is a {\em clique\/} if all
vertices in $S$ are pairwise adjacent.
The {\em independence\/} (or {\em stability\/}) {\em number\/} of $G$,
denoted by $\ga(G)$, is the largest number of vertices in an independent
set of $G$.
The {\em clique number\/} of $G$, denoted by $\omega(G)$, is the biggest
number of vertices in a clique of $G$.
The {\em complete graph\/} of order $n$, denoted by $K_n$, is a graph
of order $n$ whose vertex set is a clique. The complement of $K_n$
is the {\em empty graph\/} of order $n$.
The {\em complete bipartite graph\/} with vertex classes $V_1$ and $V_2$,
where $V_1\cap V_2=\emptyset$,
is a graph with the vertex set $V_1\cup V_2$ and the edge set consisting of
all edges $\{v_1,v_2\}$ for $v_1\in V_1$ and $v_2\in V_2$.

If $X\subseteq V(G)$, then $G[X]$ denotes the subgraph {\em induced\/}
by $G$ on $X$ (or {\em spanned\/} by $X$ in $G$
If $X,Y\subseteq V(G)$ are disjoint, then $G[X,Y]$ denotes the bipartite
graph induced by $G$ on vertex classes $X$ and $Y$, that is,
$V(G[X,Y])=X\cup Y$ and there are those edges of $G$ connecting
a vertex in $X$ with a vertex in $Y$.

We call $X\subseteq V(G)$ {\em homogeneous\/} if it is a clique
or an independent set. We call a pair of disjoint sets $X,Y\subseteq V(G)$
{\em homogeneous\/} if $G[X,Y]$ is a complete or an empty bipartite graph.

\subsubsection{General relational structures}

A {\em vocabulary\/} is a finite sequence $R_1,\ldots,R_m$ of relation
symbols along with a sequence $k_1,\ldots,k_m$ of positive integers,
where each $k_i$ is the arity of the respective $R_i$.
If $\vocab$ is a vocabulary, a finite {\em structure $G$ over $\vocab$\/}
(or {\em $\vocab$-structure $G$\/}) is a finite set $V(G)$, called the
{\em universe\/}, along with relations $R_1^G,\ldots,R_m^G$, where
$R_i^G$ has arity $k_i$. The {\em order\/} of a structure $G$ is the number
of elements in the universe $V(G)$. If $U\subseteq V(G)$, then $G$
induces on $U$ the structure $G[U]$ with universe $V(G[U])=U$
and relations $R^{G[U]}_1,\ldots,R^{G[U]}_m$ such that
$R_i^{G[U]}(\bar v)=R_i^G(\bar v)$ for every $\bar v\in U^{k_i}$.
Two $\vocab$-structures $G$ and $H$ are {\em isomorphic\/} if there is
a one-to-one map $\phi\function{V(G)}{V(H)}$, called {\em an isomorphism
from $G$ to $H$\/}, such that
$R_i^G(\bar v)=R_i^H(\phi(\bar v))$ for every $i\le m$ and all
$\bar v\in V(G)^{k_i}$. If $G$ and $H$ are isomorphic, we write $G\cong H$.
An {\em automorphism of $G$\/} is an isomorphism from $G$ to itself.
If $U\subseteq V(G)$ and $W\subseteq V(H)$,
we call a one-to-one map $\phi\function UW$ a {\em partial isomorphism
from $G$ to $H$\/} if it is an isomorphism from $G[U]$ to $H[W]$.

All these notions coincide with their standard graph-theoretic
counterparts if we consider graphs structures with a single symmetric and
anti-reflexive binary relation.

\subsection{Logic}

First order formulas are assumed to be over the set of connectives
$\{\neg,\And,\Or\}$.

\begin{definition}
A {\em sequence of quantifiers\/} is a finite word over the alphabet
$\{\exists,\forall\}$. If $S$ is a set of such sequences, then
$\exists S$ (resp.\ $\forall S$) means the set of concatenations
$\exists s$ (resp.\ $\forall s$) for all $s\in S$. If $s$ is a sequence
of quantifiers, then $\bar s$ denotes the result of replacement of all
occurrences of $\exists$ to $\forall$ and vice versa in $s$. The
 set $\bar S$
consists of all $\bar s$ for $s\in S$.

Given a first order formula $\Phi$, its set of {\em sequences of nested
quantifiers\/} is denoted by $\nest\Phi$ and defined by induction as
follows:
\begin{enumerate}
\item
$\nest\Phi=\{ \epsilon \}$ if $\Phi$ is atomic; here, $\epsilon$ denotes the
empty word.
\item
$\nest{\neg\Phi}=\overline{\nest\Phi}$.
\item
$\nest{\Phi\And\Psi}=\nest{\Phi\Or\Psi}=\nest\Phi\cup\nest\Psi$.
\item
$\nest{\exists x\Phi}=\exists\nest\Phi$ and
$\nest{\forall x\Phi}=\forall\nest\Phi$.
\end{enumerate}
\end{definition}

\begin{definition}
The {\em quantifier rank\/} of a formula $\Phi$, denoted by $\qr\Phi$
is the maximum length of a string in $\nest\Phi$.
\end{definition}

\begin{proposition}\label{prop:ld}
Let $\Phi$ be a first order formula with $\qr\Phi=k$ and suppose that none of
variables $x_1,\ldots,x_k$ occurs in $\Phi$. Then there is an equivalent
formula $\Psi$ whose bound variables are all in the set $\{x_1,\ldots,x_k\}$.
\end{proposition}

\begin{proof}
Let $m=\conn\Phi$, where $\conn\Phi$ denotes the number
of connectives in $\Phi$. We proceed by induction on $k$ and $m$.
The base cases of $k=0$, $m$ arbitrary and $m=0$, $k$ arbitrary
are straightforward. Assume the claim is true for all formulas $\Phi'$
with $\qr{\Phi'}<k$, $\conn{\Phi'}\le m$ and with
$\qr{\Phi'}\le k$, $\conn{\Phi'}<m$. If $\Phi=\neg\Phi'$, or
$\Phi=\Phi'\And\Phi''$, or $\Phi=\Phi'\Or\Phi''$, we apply the
induction assumption to $\Phi'$ and $\Phi''$, obtain equivalent formulas
$\Psi'$ and $\Psi''$, and set $\Psi=\neg\Psi'$, or $\Psi=\Psi'\And\Psi''$,
or $\Psi=\Psi'\Or\Psi''$ respectively.
If $\Phi=\exists x \Phi'$ or $\Phi=\forall x \Phi'$, we apply the
induction assumption to $\Phi'$ to obtain an equivalent formula $\Psi'$
with bound variables in $\{x_1,\ldots,x_{k-1}\}$ and rename $x$ to $x_k$.
\end{proof}

We adopt the notion of the {\em alternation number\/} of a formula
(cf.\ \cite[Definition 2.8]{Pez}).

\begin{definition}
Given a sequence of quantifiers $s$, let $\alt s$ denote the number
of occurrences of $\exists\forall$ and $\forall
\exists$ in $s$.
The {\em alternation number\/} of a first order formula $\Phi$,
denoted by $\alt\Phi$, is the maximum $\alt s$ over $s\in\nest\Phi$.
\end{definition}

\begin{definition}
Let $G$ and $G'$ be non-isomorphic structures over the same vocabulary
$\vocab$ and $\Phi$ be a first order formula over vocabulary
$\vocab\cup\{{=}\}$. We say that $\Phi$ {\em distinguishes $G$ from $G'$} if
$\Phi$ is true on $G$ but false on $G'$.
By $\D{G,G'}$ (resp.\ $\DD k{G,G'}$) we denote the minimum quantifier rank of
a formula (with alternation number at most $k$ resp.) distinguishing $G$
from $G'$. By $\V{G,G'}$ we denote the minimum $l$ such that over
the variable set $\{x_1,\ldots,x_l\}$ there is
a formula distinguishing $G$ from~$G'$.
\end{definition}

\noindent
Note that
$$
\V{G,G'}\le\D{G,G'}\le\DD k{G,G'}\le\DD{k-1}{G,G'}
$$
for every $k\ge 1$, where the first inequality follows from
Proposition~\ref{prop:ld}.

\begin{definition}
We say that $\Phi$ {\em defines\/} a structure $G$ (up to isomorphism)
if $\Phi$ distinguishes $G$ from any non-isomorphic structure $G'$
over the same vocabulary.
By $\D G$ (resp.\ $\DD kG$) we denote the minimum
quantifier rank of a formula defining $G$
(with alternation number at most $k$ resp.).
\end{definition}

\begin{definition}\label{def:var}
Let $G$ be an $\vocab$-structure. We define
$$
\V G=\max\setdef{\V{G,G'}}{G'\not\cong G},
$$
where $G'$ ranges over all $\vocab$-structures non-isomorphic with $G$.%
\footnote{In other terms, $\V G$ is equal to the minimum $l$ such that
$G$ is definable in the $l$-variable fragment of the infinitary
logic $FO_{\infty\omega}$ (see~\cite{EFl}).}
\end{definition}

\begin{proposition}\label{prop:ddd}
Let $G$ be an $\vocab$-structure.
\begin{eqnarray*}
\D G&=&\max\setdef{\D{G,G'}}{G'\not\cong G},\hspace{60mm}\mbox{}\\
\DD kG&=&\max\setdef{\DD k{G,G'}}{G'\not\cong G},
\end{eqnarray*}
where $G'$ ranges over all $\vocab$-structures non-isomorphic with $G$.
Moreover, if $G$ is a graph, one can suppose that $G'$ ranges in
the class of graphs rather than in the class of structures with
a single binary relation.
\end{proposition}

\begin{proof}
We prove the first equality; The proof of the second equality is similar.
Given $G'$ non-isomorphic with $G$, let $\Phi_{G'}$ be a formula
of minimum quantifier rank distinguishing $G$ from $G'$, that is,
$\qr{\Phi_{G'}}=\D{G,G'}$. Let $R=\max_{G'}\qr{\Phi_{G'}}$.
We have $\D G\ge R$ because $\D G\ge\D{G,G'}$ for every $G'$.
To prove the reverse inequality $\D G\le R$, notice that $G$
is defined by the formula $\Phi=\bigwedge_{G'}\Phi_{G'}$
whose quantifier rank is $R$. The only problem is that $\Phi$
is an infinite conjunction (a $FO_{\infty\omega}$-formula).
However, as it is well known, over a fixed finite vocabulary
there are only finitely many inequivalent first order formulas of
bounded quantifier rank (see e.g.\ \cite{CFI,EFl,Imm}).
We therefore can reduce $\Phi$ to a finite conjunction.

If $G$ is a graph, we do not need to consider $G'$ which are not graphs
because those are distinguished from $G$ by formula
$\exists x_1 E(x_1,x_1)\Or
\exists x_1\exists x_2(E(x_1,x_2)\And\neg E(x_2,x_1))$.
Notice that this formula has rank 2 while
no formula of rank 1 can define~$G$.
\end{proof}

\begin{proposition}\label{prop:compl}
Let $G$ and $G'$ be non-isomorphic graphs, ordinary or directed.
\begin{enumerate}
\item
$\D{G,G'}=\D{\compl{G},\compl{G'}}$,
$\DD k{G,G'}\allowbreak=\DD k{\compl{G},\compl{G'}}$, and
$\V{G,G'}=\V{\compl{G},\compl{G'}}$.
\item
$\D G=\D{\compl G}$,
$\DD kG=\DD k{\compl G}$, and
$\V G=\V{\compl G}$.
\end{enumerate}
\end{proposition}

\begin{proof}
1)
If $\Phi$ is a first order formula, let $\compl\Phi$ be the result of
putting $\neg$ in front of every occurrence of the predicate symbol $E$
in $\Phi$. As easily seen, $G\models\Phi$ iff $\compl G\models\compl\Phi$.
It remains to note that $\Phi$ and $\compl\Phi$ have the same number
of variables, the same sequences of nested quantifiers, and, consequently,
equal quantifier ranks and alternation numbers.

2) The second item follows from the first one by Proposition~\ref{prop:ddd}.
\end{proof}

\subsection{Games}\label{ss:games}

The {\em \EF\/ game\/} is played on a pair of structures of the same
vocabulary. In the case of graphs, if we want to translate the
definition below into the language of graph theory, we just have to replace
{\em structures\/} with {\em graphs}, {\em elements\/} with
{\em vertices}, and {\em universe\/} with {\em vertex set}.

\begin{definition}
Let $G$ and $G'$ be structures with disjoint universes.
The $r$-round $l$-pebble \EF\/ game on $G$ and $G'$,
denoted by $\game_r^l(G,G')$, is played by
two players, Spoiler and Duplicator, with $l$ pairwise distinct
pebbles $p_1,\ldots,p_l$, each given in duplicate. Spoiler starts the game.
A {\em round\/} consists of a move of Spoiler followed by a move of
Duplicator. At each move Spoiler takes a pebble, say $p_i$, selects one of
the structures $G$ or $G'$, and places $p_i$ on an element of this structure.
In response Duplicator should place the other copy of $p_i$ on an element
of the other structure. It is allowed to remove previously placed pebbles
to another element and place more than one pebble on the same element.

After each round of the game, for $1\le i\le l$ let $x_i$ (resp.\ $x'_i$)
denote the element of $G$ (resp.\ $G'$) occupied by $p_i$, irrespectively
of who of the players placed the pebble on this element. If $p_i$ is
off the board at this moment, $x_i$ and $x'_i$ are undefined.
If after every of $r$ rounds it is true that
$$
x_i=x_j\mbox{ iff } x'_i=x'_j\mbox{ for all }1\le i<j\le l,
$$
and the component-wise correspondence $(x_1,\ldots,x_l)$ to
$(x'_1,\ldots,x'_l)$ is a partial isomorphism from $G$ to $G'$, this is
a win for Duplicator;  Otherwise the winner is Spoiler.

The {\em $k$-alternation\/} \EF\/ game on $G$ and $G'$ is a variant
of the game in which Spoiler is allowed to switch from one structure
to another at most $k$ times during the game, i.e., in at most $k$
rounds he can choose the structure other than that in the preceding round.
\end{definition}

\noindent
The main technical tool we will use is given by the following statement.

\nopagebreak

\begin{proposition}\label{prop:loggames}
\mbox{}

\begin{enumerate}
\item
$\V{G,G'}$ equals the minimum $l$ such that Spoiler has a winning
strategy in $\game_r^l(G,G')$ for some $r$.
\item
$\D{G,G'}$ equals the minimum $r$ such that Spoiler has a winning
strategy in $\game_r^r(G,G')$.
\item
$\DD k{G,G'}$ equals the minimum $r$ such that Spoiler has a winning
strategy in the $k$-alternation $\game_r^r(G,G')$.
\end{enumerate}
\end{proposition}

\noindent
We refer the reader to \cite[Theorem 6.10]{Imm} for the proof of the first
claim, to \cite[Theorem 1.2.8]{EFl}, \cite[Theorem 6.10]{Imm},
or \cite[Theorem 2.3.1]{Spe} for the second claim, and to \cite{Pez}
for the third claim.

Proposition \ref{prop:loggames} immediately implies Propositions
\ref{prop:ld} and \ref{prop:compl} that we earlier proved syntactically.
Note that, if we prohibit removing pebbles from one vertex to another in
$\game_r^r(G,\allowbreak G')$, this will not affect the outcome of the game.

\begin{definition}
We denote the variant of $\game_r^r(G,G')$ with removing pebbles
prohibited by $\game_r(G,G')$.
\end{definition}

The examples below are obtained by simple application of
Proposition~\ref{prop:loggames}.

\begin{example}\label{ex:lower}
\mbox{}

\begin{enumerate}
\item
$\V{K_m\sqcup \compl{K_m},K_{m+1}\sqcup \compl{K_{m-1}}}=m$,
$\D{K_m\sqcup \compl{K_m},K_{m+1}\sqcup \compl{K_{m-1}}}=m+1$.
\item
$\V{K_{m+1}\sqcup \compl{K_m},K_m\sqcup \compl{K_{m+1}}}=
\D{K_{m+1}\sqcup \compl{K_m},K_m\sqcup \compl{K_{m+1}}}=m+1$.
\end{enumerate}
\end{example}

\section{Distinguishing non-isomorphic graphs}\label{s:dist}

We here prove Theorem \ref{thm:d}, that will be restated as
Theorem \ref{thm:d1}. The proof is based on the characterization
of $\D{G,G'}$ as the length of the \EF\/ game on $G$ and $G'$
given by Proposition \ref{prop:loggames}.
The most essential part of the proof is contained in Lemma \ref{lem:d1}
that gives a winning strategy for Spoiler. This lemma
will also be used in the next section to prove Theorem \ref{thm:def}.
We first introduce a couple of useful
relations between vertices of a graph that will be
intensively exploited in the course of the proofs.

\subsection{Spoiler's preliminaries}\label{ss:sprel}

\begin{definition}\label{def:sim}
We call vertices $u$ and $v$ of a graph $G$ {\em similar\/} and write
$u\sim v$ if the transposition $(uv)$ is an automorphism of $G$.
Let $[u]_G=\setdef{v\in V(G)}{v\sim u}$, $\gs_G(u)=|[u]_G|$, and
$\gs(G)=\max_{u\in V(G)}\gs_G(u)$. If the graph is clear from the context,
the subscript $G$ may be omitted. We will call the numbers $\gs(v)$
and $\gs(G)$ the {\em similarity indices\/} of the vertex $v$ and the
graph $G$ respectively.
\end{definition}
In other words, $u\sim v$ if
every third vertex $t$ is simultaneously adjacent or not to $u$ and $v$.
We will say that $t$ {\em separates\/} $u$ and $v$ if $t$ is adjacent
to exactly one of the two vertices.

\begin{lemma}\label{lem:simprop}
\mbox{}

\begin{enumerate}
\item
${\sim}$ is an equivalence relation on $V(G)$.
\item
Every equivalence class $[u]$ is a homogeneous set.
\end{enumerate}
\end{lemma}

\begin{proof}
The lemma is straightforward. The only care should be taken to
check the transitivity. Given pairwise distinct $u$, $v$, and $w$,
let us deduce from $u\sim v$ and $v\sim w$ that $u\sim w$.
For every $t\ne u,w$, we need to show that $u$ and $t$ are adjacent
iff so are $w$ and $t$. If $t\ne v$, this is true because both
adjacencies are equivalent to the adjacency of $v$ and $t$.
There remains the case that $t=v$. Then $u$ and $v$ are adjacent iff
so are $u$ and $w$ (as $v\sim w$), which in turn holds iff
so are $w$ and $v$ (as $u\sim v$).
\end{proof}

\begin{definition}
Given $X\subset V(G)$, we will denote its complement by
$\compl X=V(G)\setminus X$. Let $u,v\in\compl X$.
We write $u\xeq v$ if the identity map of
$X$ onto itself extends to an isomorphism from $G[X\cup\{u\}]$
to $G[X\cup\{v\}]$.
\end{definition}

In other words, $u\xeq v$ if these vertices have the same adjacency
pattern to $X$, i.e., $\Gg(u)\cap X=\Gg(v)\cap X$.
Clearly, $\xeq$ is an equivalence relation on $\compl X$.

\begin{definition}\label{def:calc}
$\calC(X)$ is the partition of $\compl X$ into $\xeq$-equivalence classes.
\end{definition}

Let us notice a few straightforward properties of this partition.
If $X_1\subseteq X_2$, then $\calC(X_2)$ is a refinement of
$\calC(X_1)$ on $\compl{X_2}$.
For any $X$, the $\sim$-equivalence classes restricted to $\compl X$
refine the partition $\calC(X)$.

\begin{definition}\label{def:cmax}
Let $X\subset V(G)$. We say that $X$ is {\em $\calC$-maximal\/} if
$|\calC(X\cup\{u\})|\le|\calC(X)|$ for any $u\in\compl X$.
\end{definition}

\begin{lemma}\label{lem:alpha}
Let $X\subset V(G)$ be $\calC$-maximal. Then the partition $\calC(X)$
has the following properties.
\begin{enumerate}
\item
Every $C$ in $\calC(X)$ is a homogeneous set.
\item
If $C_1$ and $C_2$ are distinct classes in $\calC(X)$ and have at least
two elements each, then the pair $C_1,C_2$ is homogeneous.
\end{enumerate}
\end{lemma}

\begin{proof}
1)
Suppose, to the contrary, that $C$ is neither a clique nor
an independent set. Then $C$ contains three vertices $u$, $v$, and $w$
such that $u$ and $v$ are adjacent but $u$ and $w$ are not.
However, if we move $u$ to $X$, then $C$ splits into two classes,
one containing $u$ and another containing $w$.
Hence the number of equivalence classes increases at least by 1,
a contradiction.

2)
Suppose that this is not true, for example, $u\in C_1$ is adjacent to
$v\in C_2$ but not to $w\in C_2$. If we move $u$ to $X$, then $C_2$
splits into two non-empty classes and $C_1\setminus\{u\}$ stays non-empty.
Again the number of equivalence classes increases, a contradiction.
\end{proof}

\begin{lemma}\label{lem:beta}
In every graph $G$, there exists a $\calC$-maximal set of vertices $X$
such that
\begin{equation}\label{eq:cmax}
|\calC(X)|\ge|X|+1.
\end{equation}
In particular,
\begin{equation}\label{eq:xsize}
|X|\le\frac{|G|-1}2.
\end{equation}
\end{lemma}

\begin{proof}
Such an $X$ can be constructed, starting from $X=\emptyset$, by repeating
the following procedure. As long as there exists $u\in\compl X$
such that $\calC(X\cup\{u\})>\calC(X)$, we move $u$ to $X$.
As soon as there is no such $u$, we arrive at $X$ which is $\calC$-maximal.
The relation \refeq{eq:cmax} is true as it holds at the beginning and is
preserved in each step. The bound \refeq{eq:xsize} follows from the
inequality $|X|+|\calC(X)|\le|G|$.
\end{proof}

\begin{definition}
Let $X\subset V(G)$.

$Y(X)$ is the union of all single-element classes in $\calC(X)$.

$Z(X)=V(G)\setminus(X\cup Y(X))$.

$\calD(X)$ is the partition of $Z(X)$ defined by
$\calD(X)=\calC(X\cup Y(X))$.
\end{definition}

\noindent
Clearly, $\calD(X)$ refines the partition induced on $Z(X)$ by $\calC(X)$.

\begin{lemma}\label{lem:cald}
If $X\subset V(G)$ is $\calC$-maximal, then every class $D$ in $\calD(X)$
consists of pairwise similar vertices. Thus, $\calD(X)$ coincides with
the partition induced on $Z(X)$ by ${\sim}$-equivalence classes.
\end{lemma}

\begin{proof}
Let $u$ and $v$ be distinct elements of the same class $D\in\calD(X)$.
These vertices cannot be separated by any vertex $t\in X\cup Y(X)$
by the definition of $\calD(X)$. Assume that they are separated
by a $t\in Z(X)$. Let $C_1$ be the class in $\calC(X)$ including $D$
and $C_2$ be the class in $\calC(X)$ containing $t$. Since $t\notin Y(X)$,
the class $C_2$ has at least one more element except $t$.
If $C_1\ne C_2$, moving $t$ to $X$ splits up $C_1$ and does not
eliminate $C_2$. If $C_1=C_2$, moving $t$ to $X$ splits up this class
and splits up or does not affect the others. In either case,
$|\calC(X)|$ increases, giving a contradiction.
\end{proof}

\begin{definition}\label{def:phieq}
Let $\phi\function X{X'}$ be a partial isomorphism from $G$ to $G'$.
Let $v\in\compl X$ and $v'\in\compl{X'}$. We call vertices $v$ and $v'$
{\em $\phi$-similar\/} and write $v\phieq v'$ if
$\phi$ extends to an isomorphism from $G[X\cup\{v\}]$ to
$G'[X'\cup\{v'\}]$.
\end{definition}
Note that, if $u\xeq v$ and $u'\xxeq v'$, then $u\phieq u'$ iff
$v\phieq v'$. This makes the following definition correct.

\begin{definition}\label{def:clphieq}
Let $\phi\function X{X'}$ be a partial isomorphism from $G$ to $G'$.
Let $C\in\calC(X)$ and $C'\in\calC(X')$. We call $C$ and $C'$
{\em $\phi$-similar\/} and write $C\phieq C'$ if
$v\phieq v'$ for some (equivalently, for all) $v\in C$ and $v'\in C'$.
\end{definition}

Notice that, if $u\phieq u'$ and $v\phieq v'$, then the relations
$u\xeq v$ and $u'\xxeq v'$ are true or false simultaneously.
It follows that the $\phi$-similarity is a matching between
the classes in $\calC(X)$ and the classes in $\calC(X')$, i.e.,
no class can have more than one $\phi$-similar counterpart in
the other graph.

\subsection{Spoiler's strategy}

\begin{definition}\label{def:add}
If $v\in V(G)$, the notation $H=G\add v$ means that
\begin{itemize}
\item
$\gs_G(v)\ge2$ and
\item
$H$ is a graph obtained from $G$ by adding a new vertex $v'$
so that $[v]_{H}=[v]_G\cup\{v'\}$.
\end{itemize}
In other words, $v'$ is similar to $v$ and adjacent to $v$ depending
on if $[v]_G$ is a clique or an independent set.
Furthermore, we define $G\add 0v=G$ and $G\add lv=(G\add (l-1)v)\add v$ for
a positive integer~$l$.
\end{definition}

\noindent
{\bf Convention.}
In the sequel, writing $H=G\add lv$ we will not require the inclusion
$V(G)\subseteq V(H)$ assuming that $H$ is an
arbitrary isomorphic copy of $G\add lv$. When considering the
\EF\/ game on $G$ and $H$, we will in addition suppose that
the vertex sets of these graphs are disjoint.

\begin{lemma}\label{lem:d1}
If $G$ and $G'$ are non-isomorphic graphs of orders $n$ and $n'$
respectively and $n\le n'$, then
\begin{equation}\label{eq:d1n4}
\DD1{G,G'}\le(n+5)/2
\end{equation}
unless $G'=G\add(n'-n)v$ for some $v\in V(G)$.
\end{lemma}

\begin{proof}
We will describe a strategy of Spoiler winning $\game_r(G,G')$ for
$r=\lfloor(n+5)/2\rfloor$ unless $G'=G\add(n'-n)v$.
The strategy splits the game in two phases.

\smallskip

\centerline{\sc Phase 1}

\smallskip

Spoiler selects a $\calC$-maximal set of vertices $X\subset V(G)$ such
that $|\calC(X)|\ge|X|+1$, whose existence is guaranteed by
Lemma \ref{lem:beta}. Denote $s=|X|$, the number of rounds in Phase 1,
and $t=|\calC(X)|$. The bounds \refeq{eq:cmax} and \refeq{eq:xsize} read
\begin{equation}\label{eq:ts}
t\ge s+1
\end{equation}
and
\begin{equation}\label{eq:ssmall}
s\le\frac{n-1}2.
\end{equation}
If Duplicator loses in Phase 1, by \refeq{eq:ssmall} this happens
within the claimed bound.

Let $X'\subseteq V(G')$ consist of the vertices selected in Phase 1
by Duplicator and $\phi\function XX'$ be the bijection defined by the
condition that $x$ and $\phi(x)$ are selected by the players in the
same round. We assume that Phase 1 finishes without Duplicator losing
and hence $\phi$ is a partial isomorphism from $G$ to $G'$.
The following useful observation is straightforward from
Definition~\ref{def:phieq}.

\begin{claim}\label{cl:phi2}
Whenever after Phase 1 Spoiler selects a vertex $v\in V(G)\cup V(G')$,
Duplicator responds with a $\phi$-similar vertex or otherwise
immediately loses.\qed
\end{claim}

\smallskip

\centerline{\sc Phase 2}

\smallskip

Denote the classes of $\calC(X)$ by $C_1,\ldots,C_t$
and the classes of $\calC(X')$ by $C'_1,\ldots,C'_{t'}$.
If there is a class $C_i$ or $C'_j$ without any $\phi$-similar
counterpart, respectively, in $\calC(X')$ or in $\calC(X)$,
then Spoiler selects a vertex in this class and wins
according to Claim \ref{cl:phi2},
making at total $s+2\le(n+3)/2$ moves and at most one alternation between
the graphs. From now on, we therefore assume that the $\phi$-similarity
determines a perfect matching between $C_1,\ldots,C_t$
and $C'_1,\ldots,C'_{t'}$, where actually $t=t'$. For the notational
convenience, we assume that $C_i\phieq C'_i$ for all $i\le t$.

Furthermore, if there is a singleton $C_i$ or $C'_j$
whose $\phi$-similar counterpart has at least two
vertices, Spoiler selects such two vertices and again wins
according to Claim \ref{cl:phi2}. We will therefore assume that
$|C_i|=1$ iff $|C'_i|=1$.
Without loss of generality, assume that $|C_i|=|C'_i|=1$ iff $i\le q$.

Denote $Y=Y(X)$ and $Y'=Y(X')$. Let $C_i=\{y_i\}$ and $C'_i=\{y'_i\}$
for $i\le q$. Thus, $Y=\{y_1,\ldots,y_q\}$ and $Y'=\{y'_1,\ldots,y'_q\}$.
Define $\phi^*\function{X\cup Y}{X'\cup Y'}$, an extension of $\phi$, by
$\phi^*(y_i)=y'_i$.

\begin{claim}\label{cl:phiiso}
$\phi^*$ is a partial isomorphism from $G$ to $G'$, unless Spoiler
wins in the next 2 moves with no alternation,
having made at total $s+2\le(n+3)/2$ moves.
\end{claim}

\begin{subproof}
For every $i\le q$, the restriction
$\phi^*\function{X\cup C_i}{X'\cup C'_i}$
is an isomorphism because $C_i$ and $C'_i$ are $\phi^*$-similar.
The restriction $\phi^*\function Y{Y'}$ should be an isomorphism as well
by the following reason.
If there are $i,j\le q$ such that $y_i$ and $y_j$
are adjacent but $y'_i$ and $y'_j$ are not or vice versa, then
Spoiler wins on the account of Claim \ref{cl:phi2} by selecting
$y_i$ and $y_j$.
\end{subproof}

We will therefore assume that $\phi^*$ is indeed a partial
isomorphism from $G$ to $G'$.
Let $Z=Z(X)$ and $Z'=Z(X')$.
Denote the classes of $\calD(X)$ by $D_1,\ldots,D_p$
and the classes of $\calD(X')$ by $D'_1,\ldots,D'_{p'}$.
Note that
\begin{equation}\label{eq:ptq2}
p\ge t-q
\end{equation}
and
$$
p=t-q\mbox{\ \ iff\ \ }\calD(X)=\{C_{q+1},\ldots,C_t\}.
$$

\begin{claim}\label{cl:phi*}
Whenever in Phase 2 Spoiler selects a vertex $v\in Z\cup Z'$,
Duplicator responds with a $\phi^*$-similar vertex or otherwise
Spoiler wins in the next round at latest,
with no alternation between $G$ and $G'$ in this round.
\end{claim}

\begin{subproof}
Let $u$ be the vertex selected by Duplicator in response to $v$
and assume that $u\notphistar v$. Suppose that $v\in Z'$
(the case of $v\in Z$ is completely similar). If $u\notin Z$,
Duplicator has already lost by Claim \ref{cl:phi2}.
If $u\in Z$, there exists a vertex
$w\in X\cup Y$ such that $u$ and $w$ are adjacent but $v$ and $\phi^*(w)$
are not or vice versa. If $w\in X$, again Duplicator has already lost.
If $w\in Y$, then in the next round Spoiler selects $\phi^*(w)$ and wins.
\end{subproof}

Claim \ref{cl:phi*} implies that every class in $\calD(X)$ or
$\calD(X')$ has a $\phi^*$-similar counterpart in, respectively,
$\calD(X')$ or $\calD(X)$ unless Spoiler wins making in Phase 2 two moves
and at most one alternation between the graphs. We will therefore
assume that this is true, that is, the $\phi^*$-similarity determines
a perfect matching between the classes $D_1,\ldots,D_p$
and $D'_1,\ldots,D'_{p'}$, where actually $p=p'$. For the notational
convenience, we assume that $D_i\phistareq D'_i$ for all $i\le p$.

\begin{claim}\label{cl:didj}
Unless Spoiler is able to win making 2 moves and at most 1 alternation
in Phase 2, the following conditions are met.
\begin{enumerate}
\item
For every $i\le p$, $D_i$ and $D'_i$ are simultaneously cliques or
independent sets.
\item
For every pair of distinct $i,j\le p$, $G[D_i,D_j]$ and $G'[D'_i,D'_j]$
are simultaneously complete or empty bipartite graphs.
\end{enumerate}
\end{claim}

\begin{subproof}
1)
Since $D_i$ consists of $X\cup Y$-similar and hence $X$-similar
vertices, by Item 1 of Lemma \ref{lem:alpha}, $D_i$ is either a clique
or an independent set. This is actually true for the class $C\in\calC(X)$
including $D_i$. If $D'_i$ has at least 2 vertices and is not a clique
or an independent set simultaneously with $D_i$, Spoiler wins in 2 moves
with 1 alternation by selecting in $D'_i$ two vertices which are
non-adjacent in the former case and adjacent in the latter case.
Indeed, if Duplicator responds with two vertices in $C$, those are
in the opposite adjacency relation. If at least one Duplicator's
response is not in $C$, he loses by Claim \ref{cl:phi2}. Note that
this argument applies, in particular, in the case that $|D'_i|\ge2$
but $|D_i|=1$.

2)
If $D_i$ and $D_j$ are included in the same class $C\in\calC(X)$,
then by Item 1 of Lemma \ref{lem:alpha}, $G[D_i,D_j]$ is either complete
or empty. Similarly to the above, complete or empty respectively
must be $G'[D'_i,D'_j]$ unless Spoiler wins in 2 moves with 1 alternation.
If $D_i$ and $D_j$ are included in different classes of $\calC(X)$,
$C^1$ and $C^2$ respectively, then, since both $C^1$ and $C^2$ have
at least 2 vertices, $G[D_i,D_j]$ is either complete or empty according
to Item 2 of Lemma \ref{lem:alpha}. If $G'[D'_i,D'_j]$ is not, respectively,
complete or empty, then Spoiler wins in 2 moves with 1 alternation by
selecting, respectively, non-adjacent or adjacent vertices, one in
$D'_i$ and another in $D'_j$. Indeed, if Duplicator responds with
one vertex in $C^1$ and another in $C^2$, those are in the opposite
adjacency relation. Otherwise Duplicator loses by Claim \ref{cl:phi2}.
\end{subproof}

Thus, in what follows we assume that the two conditions in
Claim \ref{cl:didj} are obeyed.
Together with the fact that the $D_i$'s and the $D'_i$'s are classes
of the partitions $\calC(X\cup Y)$ and $\calC(X'\cup Y')$, this implies
that each $D_i$ and each $D'_i$ consists of pairwise similar vertices
(in the sense of Definition \ref{def:sim}). Moreover, since
$D_i\phistareq D'_i$ for every $i\le p$, where $\phi^*$ is an isomorphism
from $G[X\cup Y]$ to $G'[X'\cup Y']$ and $D_i$ and $D'_i$ are
simultaneously cliques or independent sets, the graphs $G$ and $G'$ are
isomorphic iff $|D_i|=|D'_i|$ for every $i\le p$.
Since $G$ and $G'$ are supposed to be non-isomorphic, there is $D_i$
such that $|D_i|\ne|D'_i|$. We will call such a $D_i$ {\em useful\/}
(for Spoiler).

\begin{claim}\label{cl:useful}
If $D_i$ is useful and $p>t-q$, then Spoiler is able to win
having made in Phase 2 at most $\min\{|D_i|,|D'_i|\}+2$ moves
and at most 1 alternation between $G$ and $G'$. If $p=t-q$, then
$\min\{|D_i|,|D'_i|\}+1$ moves and 1 alternation suffice.
\end{claim}

\begin{subproof}
Spoiler selects $\min\{|D_i|,|D'_i|\}+1$ vertices in the larger
of the classes $D_i$ and $D'_i$. Duplicator is enforced to at least
once reply with not a $\phi^*$-similar vertex.
Then, according to Claim \ref{cl:phi*}, Spoiler wins in the next move
at latest. If $p=t-q$ and hence the $D$-classes coincide with the
$C$-classes, this extra move is not needed. This follows from
Claim \ref{cl:phi2} because in this case violation of the
$\phi^*$-similarity causes violation of the $\phi$-similarity.
\end{subproof}

Suppose that there are two useful classes, $D_i$ and $D_j$. Observe that
\begin{eqnarray}
&|D_i|+|D_j|=|Z|-\sum_{l\ne i,j}|D_l|\le (n-s-q)-(p-2)&\nonumber\\
&\le\cases{
(n-s-q)-(t-q-1)\le n-2s & if $p>t-q$,\cr
(n-s-q)-(t-q-2)\le n-2s+1 & if $p=t-q$,\cr
}&\label{eq:ptq}
\end{eqnarray}
where we use \refeq{eq:ptq2} and \refeq{eq:ts}.
It follows that one of the useful classes has at most $(n-2s+1)/2$
vertices if $p=t-q$ and at most $(n-2s)/2$ vertices if $p>t-q$.
Therefore, if $p=t-q$, Spoiler wins the game in at most
$s+(n-2s+1)/2+1=(n+3)/2$ moves and, if $p>t-q$, in at most
$s+(n-2s)/2+2=(n+4)/2$ moves, which is within the required
bound~\refeq{eq:d1n4}.

Finally, suppose that there is a unique useful class $D_m$.
According to Claim \ref{cl:useful}, Spoiler is able to win
in at most $|D_m|+2$ moves, with the total number of moves
$s+|D_m|+2$ that is within the required
bound \refeq{eq:d1n4} provided $|D_m|=1$.
Thus, we arrive at the conclusion that the bound \refeq{eq:d1n4}
may not hold true in the only case that there is exactly one
useful class $D_m$ and $|D_m|\ge 2$.
Note that we then have $n'>n$ and $|D'_m|=|D_m|+(n'-n)$.
It remains to notice that, if we remove $n'-n$ vertices
from $D'_m$, we obtain a graph isomorphic to $G$.
It follows that $G'=G\add(n'-n)v$ with $v\in D_m$.
\end{proof}

Note a direct consequence of Lemma \ref{lem:d1} and
Proposition \ref{prop:ddd}, that will be significantly
improved in the next section.

\begin{corollary}
If $\gs(G)=1$, then $\DD1G\le(n+5)/2$.
\end{corollary}

We now restate and prove Theorem \ref{thm:d} from the introduction.

\begin{theorem}\label{thm:d1}
If $G$ and $G'$ are non-isomorphic and have the same order $n$, then
\begin{equation}\label{eq:d1n3}
\DD1{G,G'}\le(n+3)/2.
\end{equation}
\end{theorem}

\begin{proof}
Lemma \ref{lem:d1} immediately gives an upper bound of $(n+5)/2$,
which is a bit worse than we now claim.
To improve it, we go trough lines of the proof of Lemma \ref{lem:d1}
but make use of the equality $n=n'$. The latter causes the following
changes. Since $n'=n$, there must be at least two useful classes,
$D_i$ and $D_j$, such that $|D_i|<|D'_i|$ and $|D_j|>|D'_j|$.
If $p=t-q$, the bound of $(n+3)/2$ has been actually proved, and
we only need to tackle the case that $p>t-q$. Similarly to \refeq{eq:ptq},
we have
$$
2|D_i|+2|D'_j|+2\le|D_i|+|D_j|+|D'_i|+|D'_j|\le2\of{(n-s-q)-(p-2)}\le2(n-2s).
$$
It follows that at least one of $|D_i|$ and $|D'_j|$ does not
exceed $(n-2s-1)/2$. By Claim \ref{cl:useful}, Spoiler wins in totally
at most $s+(n-2s-1)/2+2=(n+3)/2$ moves.
\end{proof}

In the conclusion of this section, we state a lemma for further use
in Section \ref{s:bdeg}. This lemma is actually a corollary from the proof
of Lemma \ref{lem:d1}. More precisely, it is a variant of
Claim \ref{cl:useful}, where we take into account Lemma~\ref{lem:cald}.

\begin{lemma}\label{lem:usefulgs}
Let $G$ and $G'$ be arbitrary non-isomorphic graphs.
Suppose that $X\subset V(G)$ is $\calC$-maximal. Then
Spoiler wins the 1-alternation \EF\/ game on $G$ and $G'$
in at most $|X|+\max_{v\notin X}\gs_G(v)+2$ rounds.\noproof
\end{lemma}

\section{Defining a graph}\label{s:def}

Our next goal is to prove Theorem \ref{thm:def}. Below this theorem is
restated as Theorem \ref{thm:def1} after precisely defining the class
of graphs whose members have a larger but efficiently computable $\D G$
(see Definition \ref{def:cls}).
The proof is based on the following four lemmas.

\begin{lemma}\label{lem:sim}
Let $x_i$ (resp.\ $x'_i$) denote the vertex of $G$ (resp.\ $G'$) selected
in the $i$-th round of\/ $\game_r(G,G')$.
Then, as soon as a move of Duplicator
violates the condition that $x_i\sim x_j$ iff $x'_i\sim x'_j$, Spoiler
wins either immediately or in the next move possibly with one
alternation between the graphs.
\end{lemma}

\begin{proof}
Suppose, for example, that Duplicator selects $x'_j$ so that
$x'_i\not\sim x'_j$ while $x_i\sim x_j$ for some $i<j$.
Suppose that the correspondence between
the $x_m$'s and the $x'_m$'s, $1\le m\le j$, is still a partial isomorphism.
Then there is $y\in V(G')$ adjacent to exactly one of $x'_i$ and $x'_j$.
Note that such $y$ could not be selected by the players previously.
In the next move Spoiler selects $y$ and wins, whatever the move of
Duplicator is.
\end{proof}

\begin{lemma}\label{lem:close}
Let $|G|=n$, $v\in V(G)$, $\gs_G(v)=s$, and $G'=G\add lv$ with $l\ge1$. Then
$$
s+1\le\V{G,G'}\le\DD1{G,G'}\le s+1+\frac{n+1}{s+1}\le\cases{
(n+5)/2\mbox{\ \ for\ \ }1\le s\le(n-1)/2,&\cr
s+3-1/(n/2+1)\mbox{\ \ for\ \ }s\ge n/2.&\cr
}
$$
\end{lemma}

\begin{proof}
The lower bound is given by the following strategy for Duplicator
in $\game_r^s(G,G')$. Whenever Spoiler selects a vertex outside $[v]$
in either graph, Duplicator selects its copy in the other graph.
If Spoiler selects an unoccupied vertex similar to $v$,
then Duplicator selects an arbitrary unoccupied vertex similar to $v$
in the other graph. Clearly, this strategy preserves the isomorphism
arbitrarily long, that is, is winning for every $r$.

The upper bound for $\DD1{G,G'}$ is ensured by the following Spoiler's
strategy winning in the 1-alternation $\game_r(G,G')$ for
$r=\lfloor s+1+\frac{n+1}{s+1}\rfloor$.
In the first round Spoiler selects a vertex in $[v]_{G'}$.
Suppose that Duplicator replies with a vertex in $[u]_G$.

\case 1 {$|[u]_G|\le s$}
Spoiler continues to select vertices in $[v]_{G'}$.
In the $(s+1)$-th round at latest, Duplicator selects a vertex
outside $[u]_G$. Spoiler wins in the next move by Lemma \ref{lem:sim},
having made at most $s+2$ moves and one alternation.

\case 2 {$|[u]_G|\ge s+1$}
Spoiler selects one vertex in each similarity class of $G'$
containing at least $s+1$ vertices. Besides $[v]_{G'}$, there
can be at most $\frac{n-s}{s+1}$ such classes. At latest in the
$\lfloor\frac{n-s}{s+1}+1\rfloor$-th round Duplicator
selects either another vertex in a class with an already selected
vertex (then Spoiler wins in one extra move by Lemma \ref{lem:sim})
or a vertex in $[w]_G$ with $|[w]_G|\le s$. In the latter case
Spoiler selects $s$ more vertices in the corresponding class of $G'$.
Duplicator is forced to move outside $[w]_G$ and loses in the next
move by Lemma \ref{lem:sim}. Altogether there are made at most
$\lfloor\frac{n-s}{s+1}+1\rfloor+s+1\le s+1+\frac{n+1}{s+1}$ moves.

If $s\ge n/2$, the last inequality of the lemma is straightforward and,
if $2\le s\le(n-1)/2$, it follows from
the fact that the function $f(x)=x+\frac{n+1}x$ on $[2,(n+1)/2]$
attains its maximum at the endpoints of this range.
\end{proof}

\noindent
Using Lemma \ref{lem:close}, Lemma \ref{lem:d1} can now be refined.

\begin{lemma}\label{lem:upper}
If $G$ and $G'$ are graphs of orders $n\le n'$, then
$$
\DD1{G,G'}\le(n+5)/2
$$
unless
\begin{equation}\label{eq:s}
\gs(G)\ge n/2\mbox{ and }
G'=G\add(n'-n)v\mbox{ for some }v\in V(G)\mbox{ with }\gs_G(v)=\gs(G).
\end{equation}
In the latter case we have
\begin{equation}\label{eq:s12}
\gs(G)+1\le\V{G,G'}\le\DD1{G,G'}\le \gs(G)+2.
\end{equation}
\end{lemma}

\noindent
Note that the condition \refeq{eq:s}
determines $G'$ up to isomorphism with two exceptions if $n$ is even.
Namely, for $G=K_m\sqcup \compl{K_m}$ and $G=\compl{K_m\sqcup \compl{K_m}}$
there are two ways to extend $G$ to~$G'$.

The gap between the bounds \refeq{eq:s12} can be completely closed.

\begin{lemma}\label{lem:exact}
Let $G$ and $G'$ be graphs of orders $n\le n'$.
Assume the condition \refeq{eq:s}. Then
$
\V{G,G'}=\gs(G)+1
$
for all such $G$ and $G'$,
$
\DD0{G,G'}=\gs(G)+1
$
if $[v]_G$ is an inclusion maximal homogeneous set, and
$
\D{G,G'}=\gs(G)+2
$
if $[v]_G$ is not.
\end{lemma}

\noindent
Using the inequality $|[v]_G|\ge n/2$ and the mutual similarity
of vertices in $[v]_G$, one can easily show that $[v]_G$ is a
maximal (with respect to the inclusion) clique or independent set
iff $|[v]_G|=\ga(G)$
or $|[v]_G|=\omega(G)$ respectively. However, the former condition
is more preferable because it is efficiently verifiable.

\begin{proof}
For simplicity we assume that $[v]_G$ is an independent set.
Otherwise we can switch to $\compl G$ and $\compl{G'}$ by
Proposition \ref{prop:compl}. Denote $s=\gs(G)=|[v]_G|$.

\case 1 {$[v]_G$ is maximal independent}
We show the bound $\DD0{G,G'}\le s+1$ by describing Spoiler's strategy
winning the 0-alternation $\game_{s+1}(G,G')$. Spoiler selects $s+1$
vertices in $[v]_{G'}$.
Duplicator is forced to select at least one vertex $u_1\in[v]_G$
and at least one vertex $u_2\notin[v]_G$. Since $[v]_G$ is a maximal
independent set, $u_1$ and $u_2$ are
adjacent and this is Spoiler's win.

\case 2 {$[v]_G$ is not maximal}
We first show the bound $\V{G,G'}\le s+1$ by describing Spoiler's
strategy winning $\game_{s+2}^{s+1}(G,G')$. As in the preceding case,
Spoiler selects $s+1$ vertices in $[v]_{G'}$ and there are $u_1\in[v]_G$
and $u_2\notin[v]_G$ selected in response by Duplicator.
Assume that $u_1$ and $u_2$ are not adjacent for otherwise
Duplicator loses immediately. Since $u_1$ and $u_2$ are not similar,
there is $u\in V(G)\setminus\{u_1,u_2\}$ adjacent to exactly one of
$u_1$ and $u_2$. It follows that $u\notin[v]_G$. Note that $u$ could not be
selected by Duplicator in the first $s+1$ rounds without immediately
losing. Therefore, Duplicator has selected in $[v]_G$ at least two
vertices, say, $u_0$ and $u_1$. In the $(s+2)$-th round
Spoiler removes the pebble from $u_0$ to $u$ and
wins because the counterparts of $u_1$ and $u_2$
in $G'$ are similar and hence equally adjacent or non-adjacent
to any counterpart of~$u$.

We now show the bound $\D{G,G'}>s+1$ by describing Duplicator's
strategy winning $\game_{s+1}(G,G')$. Whenever Spoiler selects
a vertex of either graph, Duplicator selects its copy in the
other graph, with the convention that the copy of a vertex in $[v]_{G'}$
is an arbitrary unselected vertex in $[v]_G$. This is impossible
in the only case when Spoiler selects $s+1$ vertices all in $[v]_{G'}$.
Then Duplicator, in addition to $s$ vertices of $[v]_G$, selects one more
vertex extending $[v]_G$ to a larger independent set.
\end{proof}

\begin{definition}\label{def:cls}
$\cls$ is the class of graphs $G$ with $\gs(G)>(|G|+3)/2$.
$\cls_1$ is the class of graphs $G$ with $\gs(G)>(|G|+3)/2$ such that
the largest similarity class is an inclusion maximal homogeneous set.
$\cls_2$ is the class of graphs $G$ with $\gs(G)>(|G|+1)/2$ such that
the largest similarity class is not an inclusion maximal homogeneous set.
\end{definition}

\begin{theorem}\label{thm:def1}
$\V{G}\le(|G|+5)/2$ with the exception of all graphs in $\cls$.
If $G\in\cls$, then $\V G=\gs(G)+1$.

$\DD1{G}\le(|G|+5)/2$ with the exception of all graphs in $\cls_1\cup\cls_2$.
If $G\in\cls_1$, then $\D G=\gs(G)+1$;
If $G\in\cls_2$, then $\D G=\gs(G)+2$.
\end{theorem}

\begin{proof}
We prove the theorem for $\V G$; the proof for $\D G$ is completely
similar. Recall that $\V G=\max\setdef{\V{G,G'}}{G'\not\cong G}$.
We consider two cases.

\case 1 {$\gs(G)<|G|/2$}
For every $G'\not\cong G$ we have $\V{G,G'}\le(|G|+5)/2$
by Lemma \ref{lem:upper}. Since $G\notin\cls$, the theorem in Case 1
is true.

\case 2 {$\gs(G)\ge|G|/2$}
If $G'=G\add lv$ for some $l\ge 1$ and $v\in V(G)$ with
$\gs_G(v)=\gs(G)$, then $\V{G,G'}=\gs(G)+1$ by Lemma \ref{lem:exact}.
By the definition of $\cls$, we therefore have $\V{G,G'}\le(|G|+5)/2$
if $G\notin\cls$ and $\V{G,G'}>(|G|+5)/2$ if $G\in\cls$.

If $G=G'\add lv$ for some $l\ge 1$, $G'$ with $\gs(G')\ge|G'|/2$,
and $v\in V(G')$ with $\gs_{G'}(v)=\gs(G')$, then $\V{G,G'}=\gs(G')+1$
by Lemma \ref{lem:exact} and hence $\V{G,G'}<\gs(G)+1$.

If $G'$ is any other graph non-isomorphic with $G$, then
$\V{G,G'}\le(|G|+5)/2$ by Lemma \ref{lem:upper}.

Summarizing, if $G\notin\cls$, we have $\max_{G'}\V{G,G'}\le(|G|+5)/2$
and, if $G\in\cls$, we have $\max_{G'}\V{G,G'}=\gs(G)+1>(|G|+5)/2$.
Thus, in Case 2 the theorem is also true.
\end{proof}

A variant of Theorem \ref{thm:def1} was stated in the introduction
as Theorem \ref{thm:def}. To link the two theorems, in the latter
we should set $\cl=\cls_1\cup\cls_2$. The efficiency statements
of Theorem \ref{thm:def} are due to the following lemma.
Referring to efficient algorithms, we mean
random access machines whose running time on graphs of order $n$, represented
by adjacency matrices, is~$O(n^2\log n)$.

\begin{lemma}\label{lemma:eff}
\mbox{}

\begin{enumerate}
\item
There is an efficient algorithm that, given $G$, finds the partition
of $V(G)$ into classes of pairwise similar vertices.
\item
Given $G$, the number $\gs(G)$ is efficiently computable.
\item
The classes $\cls$, $\cls_1$, and $\cls_2$ defined in Definition
\ref{def:cls} are efficiently recognizable.
\end{enumerate}
\end{lemma}

\begin{proof}
Notice that non-adjacent vertices are in the same similarity class iff
the corresponding rows of the adjacency matrix are identical. Thus,
in order to find similarity classes containing more than one element,
it suffices, using the standard $O(n\log n)$-comparison sorting,
to arrange rows of the adjacency matrices of $G$ and $\compl G$
in the lexicographic order. This proves Item 1. The other two
are its direct consequences.
\end{proof}

\begin{remark}
An analysis of the proofs shows that Theorem \ref{thm:def1} is
even more constructive: Given a graph $G$, one can
efficiently construct its defining formula
whose quantifier rank is as small as possible if $G\in\cls_1\cup\cls_2$ and
does not exceed $(n+5)/2$ if $G\notin\cls_1\cup\cls_2$.
\end{remark}

\begin{remark}\label{rem:arbstr}
Note that Definition \ref{def:sim} of the similarity relation makes
sense for an arbitrary structure.
Lemma \ref{lem:sim} generalizes over any class of $\vocab$-structures,
for an arbitrary vocabulary $\vocab$: If precisely one of the relations
$x_i\sim x_j$ and $x'_i\sim x'_j$ holds, then Spoiler is able to win
in at most $k-1$ next moves, where $k$ is the maximum relation arity
of $\vocab$. If, say, $x'_i\not\sim x'_j$ but $x_i\sim x_j$, the only what
Spoiler needs to do is to exhibit $u_1,\ldots,u_{k-1}\in V(G')$
such that there is a sequence consisting of elements $u_1,\ldots,u_{k-1}$
and a variable $x$ that satisfies some relation with $x=x_i$
and does not satisfy the same relation with $x=x_j$.

Lemma \ref{lem:close}, whose proof uses only Lemma \ref{lem:sim}
and the definition of similarity classes, carries over to structures
with maximum relation arity $k$ giving bounds
$$
s+1\le\V{G,G'}\le\DD1{G,G'}\le s+k-1+\frac{n+1}{s+1}.
$$
\end{remark}

\section{Distinguishing graphs by zero-alternation formulas}\label{s:d0}

Theorem \ref{thm:d1} is proved in a stronger form: The
class of distinguishing formulas is restricted to those
with alternation number 1. We now further restrict the
alternation number to the smallest possible value of 0.
In terms of the \EF\/ game, we restrict the ability of Spoiler
to alternate between graphs during play (see Item 3 of
Proposition \ref{prop:loggames}). Moreover, let us call a formula
{\em existential\/} (resp.\ {\em universal}) if it is in the negation
normal form and all quantifiers in it are existential (resp.\ universal).
It is easy to prove that, if a graph $G$ is distinguished from
another graph $G'$ by a formula with alternation number 0, then
$G$ is distinguished from $G'$ by either existential or universal
formula of the same quantifier rank. Somewhat surprizingly, this
restriction of the class of distinguishing formulas turns out
not so essential in the worst case.

\begin{theorem}\label{thm:D0}
If $G$ and $G'$ are non-isomorphic and have the same order $n$, then
$\DD0{G,G'}\le(n+5)/2$.
\end{theorem}

\begin{proof}
We will describe a strategy for Spoiler winning the 0-alternation
game $\game_{\lfloor(n+5)/2\rfloor}(G,G')$.
Given a set of vertices $X$ in a graph and a partial isomorphism
$\phi\function X{X'}$ to another graph, we will use the notions
introduces in Section \ref{ss:sprel}: the partitions $\calC(X)$
and $\calD(X)$, the set $Y(X)$, and the $\phi$-similarity relation
$\phieq$. We set the following notation:
\begin{eqnarray*}
s(X)&=&|X|,\\
t(X)&=&|\calC(X)|,\\
c(X)&=&\max\setdef{|C|}{C\in\calC(X)},\\
d(X)&=&\max\setdef{|D|}{D\in\calD(X)}.
\end{eqnarray*}
For brevity, we will not indicate the dependence on $X$,
writing merely $Y$, $s$, $t$, $c$, and $d$.

At the start of the game Spoiler, over all choices of $H=G$ or
$H=G'$, and of $X\subset V(H)$ with $t\ge s+1$ takes one which
\begin{description}
\item[\rm(criterion 1)]
first maximizes $s$;
\item[\rm(criterion 2)]
then, if there is still some choice, minimizes $c$;
\item[\rm(criterion 3)]
finally, minimizes $d$.
\end{description}
Let us assume $H=G$. As $s+t\le n$, we have $s\le(n-1)/2$.

Spoiler selects all vertices in $X$ in any order.
Denote the set of vertices of $G'$ selected in response by Duplicator
by $X'$. Let $t'=t(X')$ and $c'=c(X')$. Assume that Duplicator has not
lost up to now, that is, has managed to maintain the partial isomorphism
$\phi\function X{X'}$. Let $C_1,\ldots,C_t$ (resp.\ $C'_1,\ldots,C'_{t'}$)
be all classes in $\calC(X)$ (resp.\ $\calC(X')$).

If there is a class $C_i$ without $\phi$-similar counterpart in $\calC(X')$,
Spoiler wins in one move by selecting a vertex in $C_i$, having made
$s+1\le(n+1)/2$ moves at total. We therefore suppose that $t'\ge t$
and, for every $i\le t$, the classes $C_i$ and $C'_i$ are $\phi$-similar.
Thus, $t'\ge s+1=s'+1$ and, by Criterion 2 of the choice of $(H,X)$,
we conclude that there is $C'_m$ such that
$|C_i|\le|C'_m|$ for all $i$.

If $t'>t$, define $C_i=\emptyset$ for all $t<i\le t'$.
Suppose first that for some $i\le t'$ we have
$|C_i|\not=|C'_i|$. As $G$ and $G'$ have the same order, there must be
an $i\le t'$ with $|C_i|> |C'_i|$. Spoiler wins by selecting $|C'_i|+1$
vertices inside $C_i$. Observe that $|C'_i|<|C_i|\le |C'_m|$ and that
$$
2\,|C'_i|+1\le |C'_i|+|C'_m|\le (n-s')-(t'-2)\le n-2s+1.
$$
The total number of Spoiler's moves is therefore at most
$s+|C'_i|+1\le n/2+1$, within the required bound.

Suppose from now on, that $t'=t$ and for any $i\le t$ we have $|C_i|=|C'_i|$.
In particular,
\begin{equation}\label{eq:cc}
c=c'.
\end{equation}
Without loss of generality, assume that $|C_i|=|C'_i|=1$ precisely
for $i\le q$. Note that $Y=\bigcup^q_{i=1}C_i$ and $Y'=\bigcup^q_{i=1}C'_i$,
where $Y'=Y(X')$. Similarly to the proof of Lemma \ref{lem:d1}, we
extend $\phi$ to $\phi^*\function{X\cup Y}{X'\cup Y'}$ by the condition
that $\phi^*$ maps each $C_i$ with $i\le q$ onto $C'_i$.
Similarly to Claim \ref{cl:phiiso}, if $\phi^*$ is not an isomorphism
from $G[X\cup Y]$ to $G'[X'\cup Y']$, then Spoiler wins by selecting 2
vertices in $Y$, having made altogether $s+2\le(n+3)/2$ moves.
In the sequel we therefore suppose that $\phi^*$ is a partial
isomorphism from $G$ to $G'$. We will make use of the following
observation, provable similarly to Claim \ref{cl:phi*}.
Let $Z=V(G)\setminus(X\cup Y)$.

\begin{claim}\label{cl:cl}
{}From now on, whenever Spoiler selects a vertex $v\in Z$,
Duplicator responds with a $\phi^*$-similar vertex or otherwise
loses in the next round at latest with no alternation.\qed
\end{claim}

Let $D_1,\ldots,D_p$ (resp.\ $D'_1,\ldots,D'_{p'}$)
be all classes in $\calD(X)$ (resp.\ $\calD(X')$).
We now claim that every class $D_i$ has a $\phi^*$-similar
counterpart in $\calD(X')$ or otherwise Spoiler wins in at most 2 next
moves with no alternation, having made altogether at most
$s+2\le(n+3)/2$ moves. Indeed, if a $D_i$ has no $\phi^*$-similar
counterpart, Spoiler selects a vertex in the $D_i$ and wins either
immediately or in the next move by Claim \ref{cl:cl}. We hence will assume that
$p'\ge p$ and, for all $i\le p$, the classes $D_i$ and $D'_i$ are
$\phi^*$-similar. If $p'>p$, define $D_i=\emptyset$ for all
$p<i\le p'$.

We now show that each class in $\calD(X)$ or $\calD(X')$
consists of pairwise similar vertices as defined in
Definition \ref{def:sim}. Suppose, to the contrary, that
vertices $u$ and $v$ lie in the same $D_i$ and some $w$ is connected
to one of $u,\,v$ but not to the other. By
the definition of $D_i$ such $w$ must lie in $Z$; but then moving $w$
to $X$ we increase $t$ at least by one.
Indeed, if $w$ belongs to the same $D_i$, the class $C^1\in\calC(X)$
including $D_i$ splits up into at least two subclasses, containing
$u$ and $v$ respectively, while no class in $\calC(X)$ disappears.
If $w$ belongs to another $D_j$, the class $C^1$ splits up as well,
while the class $C^2\in\calC(X)$ including $D_j$ still stays because
it has at least two elements.
Since the relation $t\ge s+1$ is preserved, we get a contradiction
with Criterion 1 in the choice of $(H,X)$. The same argument applies
for $\calD(X')$.

It follows that, for any distinct $i,j\le p$, each of $G[D_i]$,
$G'[D'_i]$, $G[D_i,D_j]$ and $G'[D'_i,D'_j]$ is either complete or
empty. The same is true about every $G[D_i,\{v\}]$ and $G'[D'_i,\{v'\}]$
for $v\in X\cup Y$ and $v'\in X'\cup Y'$.
We now claim that, for every $i\le p$, $j\le p$ such that $j\ne i$,
and $v\in X\cup Y$,
\begin{enumerate}
\item
$G[D_i]$ with at least 2 vertices is complete iff $G'[D'_i]$ is,
\item
$G[D_i,D_j]$ is complete iff $G'[D'_i,D'_j]$ is, and
\item
$G[D_i,\{v\}]$ is complete iff $G'[D'_i,\{\phi^*(v)\}]$ is
\end{enumerate}
or otherwise Spoiler wins in at most 3 next moves with no alternation,
having made altogether $s+3\le(n+5)/2$ moves.
For example, consider the case that $G[D_i]$ has at least 2 vertices
and is complete but $G'[D'_i]$ is empty. Then Spoiler selects two
vertices in $D_i$. If both Duplicator's responses are in $D'_i$,
he loses immediately. Otherwise Duplicator responds at least once
with a vertex which is not $\phi^*$-similar. Then Spoiler wins in the next
move according to Claim~\ref{cl:cl}.

We therefore suppose that the above three conditions are obeyed
for all $i,j\le p$ and $v\in X\cup Y$. It follows that,
if $|D_i|=|D'_i|$ for all $i\le p'$
and, in particular, $p'=p$, then $G$ and $G'$ should be isomorphic.
Since this is not so, there is $l\le p'$ such that $|D_l|\ne|D'_l|$.
As $G$ and $G'$ have the same order, we can assume that
\begin{equation}\label{eq:1op}
|D_l|>|D'_l|.
\end{equation}
Note that $p'>t-q$ for else the $D'$-classes are identical
to the $C'$-classes, which contradicts \refeq{eq:1op}. Thus
\begin{equation}\label{eq:psq}
p'\ge s+2-q.
\end{equation}

It follows from \refeq{eq:cc} and Criterion 3 of the choice of $(H,X)$
that there exists $k\le p'$ such that $|D_i|\le|D'_k|$ for all $i$.
We have $|D'_l|<|D_l|\le |D'_k|$, so
$$
2\,|D'_l| +1 \le |D'_l|+|D'_k|\le (n-s-q)-(p'-2)\le n-2s,
$$
where the latter inequality follows from~\refeq{eq:psq}.

Now, Spoiler selects $|D'_l|+1$ vertices inside $D_l$. Duplicator
cannot reply to this with all moves in $D'_l$ and hence replies at least
once with a vertex which is not $\phi^*$-similar.
According to Claim \ref{cl:cl},
Spoiler wins either immediately or in the next round.
The total number of moves is at most
$$
s+|D'_l|+1+1\le s+\frac{n-2s-1}2+2=\frac{n+3}2,
$$
as required.
\end{proof}

\section{Defining graphs of bounded degree}\label{s:bdeg}

The {\em degree\/} of a vertex $v$ in a graph, denoted by $\deg(v)$,
is the number of edges incident to $v$. The {\em maximum degree\/}
of a graph $G$ is defined by $\Delta(G)=\max_{v\in V(G)}\deg(v)$.
The {\em distance\/} between vertices $v$ and $u$ in a graph, $\dist(v,u)$,
is the smallest number of edges in a path from $v$ to $u$.
If $U\subseteq V(G)$, then $\dist(v,U)=\min_{u\in U}\dist(v,u)$.
Recall that the similarity index $\gs_G(v)$ of a vertex $v$ is
defined in Definition~\ref{def:sim}.

\begin{lemma}\label{lem:gsvsdelta}
If $v$ is a non-isolated vertex of a graph $G$, then
\begin{equation}\label{eq:gsgv}
\gs_G(v)\le\Delta(G)+1.
\end{equation}
\end{lemma}

\begin{proof}
By Lemma \ref{lem:simprop}, the similarity class $[v]_G$ is
either a clique or an independent set. If it is a clique, then the
bound \refeq{eq:gsgv} is clear. Otherwise there must exist a vertex
$u\notin[v]_G$ adjacent to $v$. As $u$ is adjacent to every vertex
in $[v]_G$, we have $\gs_G(v)\le\deg(u)\le\Delta(G)$ in this case.
\end{proof}

Recall that, while a formula defining a graph $G$ distinguishes $G$
from {\em all\/} non-isomorphic graphs, a formula
identifying $G$ distinguishes $G$ from all non-isomorphic graphs
of {\em the same order}. While the minimum quantifier rank of a defining
formula with alternation number at most $k$ is denoted by $\DD kG$,
for an identifying formula it is denoted by $\II kG$.
By Proposition \ref{prop:ddd} we have
$$
\DD kG=\max\setdef{\DD k{G,G'}}{G'\not\cong G},
$$
and likewise
$$
\II kG=\max\setdef{\DD k{G,G'}}{G'\not\cong G,\ |G'|=|G|}.
$$

\begin{theorem}\label{thm:bdeg}
Let $d\ge2$.
If $G$ is a graph of order $n$ with $\Delta(G)=d$ that has no
isolated vertex and no isolated edge, then
$$
\DD1G\le c_d n+d^2+d+7/2
$$
for a constant $c_d=\frac12-\frac1{10}d^{-2d-5}$.
If $G$ is an arbitrary graph of order $n$ with $\Delta(G)=d$, then
the same bound holds for $\II1G$.
\end{theorem}
The constant $c_d$ as stated in the theorem is far from being
best possible. We do not try to optimize it; Our goal is more
moderate, just to show the existence of a $c_d$ strictly less
than $1/2$.
A tight, up to an additive constant, bound is easy to find for $d=2$.
If $\Delta(G)=2$ and $G$ has no isolated vertices and edges,
the graph is a sum of paths and cycles.
Using bounds $\DD0{C_n,C_m}\le\log_2n+O(1)$ and
$\DD0{P_n,P_m}\le\log_2n+O(1)$ for $m\ne n$
(e.g.\ \cite[Theorem 2.1.2]{Spe}),
one can show that $\DD1G\le n/3+O(1)$.

\begin{proofof}{Theorem \ref{thm:bdeg}}
We will prove the bounds for $\DD1G$ and $\II1G$ in parallel.
We actually have to estimate $\DD1{G,G'}$ in two cases:
\begin{enumerate}
\item
$G$ has no isolated vertices and edges; $G'\not\cong G$ has arbitrary order.
\item
$G$ is arbitrary; $G'\not\cong G$ has order $n$.
\end{enumerate}
In fact, almost all proof will go through for the most general case
that $G'\not\cong G$, with no other assumptions. Only once we will need
the condition $|G'|=|G|$ for $G$ with isolated vertices or edges,
and this will be explicitly stated.

Referring to Item 3 of Proposition \ref{prop:loggames},
we design a strategy for Spoiler winning the 1-alternation
\EF\/ game on $G$ and $G'$ in at most $c_d n+d^2+d+7/2$ rounds.
Clearly, we may assume that $\Delta(G')\le d$ for otherwise
Spoiler wins in at most $d+2$ moves by selecting a star
$K_{1,d+1}$ in $G'$.

A {\em component\/} of a graph is a maximal connected induced subgraph.
We call a component {\em small\/} if it has at most $d^2+1$
vertices. Throughout the proof we use the following notation.

\hangindent=2\parindent \hangafter=0
\noindent
$A\subseteq V(G)$ consists of all vertices in small components.\\
$B_1\subseteq V(G)$ consists of all isolated vertices.\\
$B_2\subseteq V(G)$ consists of all vertices in isolated edges.\\
$B=B_1\cup B_2$.\\
$a=|A|$.\\
$b_1=|B_1|$.\\
$b_2=|B_2|/2$.\\
$A'$, $B'_1$, $B'_2$, $a'$, $b'_1$, and $b'_2$ are similarly
defined for $G'$.

\noindent
We set $\tau=\frac15d^{-2d-5}$.
Spoiler will choose one of two strategies depending on how large
or small $a$ is.

\medskip

\centerline{{\sc Strategy 1} (applicable if $a\ge\tau n$)}

\smallskip

\case 1{$G[A]\cong G'[A']$}

\noindent
Spoiler plays outside $A$ and $A'$ using the strategy for the game
on non-isomorphic graphs $G[\compl A]$ and $G'[\compl{A'}]$
described in the proof of Lemma \ref{lem:d1}.
If Duplicator never moves in $A$ or $A'$, Spoiler wins, according to
Lemmas \ref{lem:beta}, \ref{lem:usefulgs}, and \ref{lem:gsvsdelta},
in at most $(n-a-1)/2+(d+1)+2=(n-a)/2+d+5/2$ moves.
If Duplicator makes a move in $A$ or $A'$, Spoiler, who has selected in
this round a vertex $v$ in a component with more than $d^2+1$ vertices,
wins by selecting a set of $d^2+2$ vertices that includes $v$
and spans a connected subgraph. Thus, Spoiler needs at most
$\frac{n-a}2+d^2+d+\frac72\le (\frac12-\frac12\tau)n+d^2+d+\frac72$
moves to win.

\smallskip

\case 2{$G[A\setminus B]\not\cong G'[A'\setminus B']$}

\noindent
Spoiler enforces play in $A\setminus B$ and $A'\setminus B'$.
He starts in $G$ if $G[A\setminus B]$ has at least as many components
as $G'[A'\setminus B']$ has and in $G'$ otherwise.
Without loss of generality assume the former.

Spoiler selects one vertex in each component of $G[A\setminus B]$.
This takes at most $n/3$ moves as every component of $G[A\setminus B]$
has at least 3 vertices. Spoiler keeps doing so until one of the
following happens.
\begin{enumerate}
\item
Duplicator moves in $B'$. Then Spoiler wins in at most 2 extra moves.
\item
Duplicator moves outside $A'$. Then Spoiler switches to $G'$ and wins
in at most $d^2+1$ extra moves by selecting a connected subgraph
spanned by $d^2+2$ vertices.
\item
While Spoiler selects a vertex in a component $C$ of $G[A\setminus B]$,
Duplicator responds with a vertex in a component $C'$
of $G'[A'\setminus B']$ such that $C'\not\cong C$. Then Spoiler
wins in at most $d^2$ extra moves by selecting all vertices of $C$
if $|C|\ge|C'|$ or all vertices of $C'$ otherwise.
\end{enumerate}
It is clear that one of the three situations must happen sooner or later.
Thus, Spoiler wins in at most $n/3+d^2+1$ moves with at most 1 alternation
between $G$ and $G'$.

\smallskip

If neither Case 1 nor Case 2 takes place, then
$$
G[B]\not\cong G'[B'].
$$

\case 3{$b_1\ne b'_1$ and $b_2\ne b'_2$}

\noindent
Since $b_1+2b_2\le n$, we have $b_1\le n/3$ or $b_2\le n/3$.
If $b_1\le n/3$ , Spoiler selects $\min\{b_1,b'_1\}+1$ isolated vertices
in the graph containing larger number of them and wins in the next move
with alternation between the graphs.
If $b_2\le n/3$, Spoiler selects one vertex in each of
$\min\{b_2,b'_2\}+1$ isolated edges in one of the graphs, where this
is possible, and then wins in at most 2 next moves with possibly
1 alternation. At total, at most $n/3+3$ moves are needed.

\smallskip

It remains to tackle the situation when exactly one of the inequalities
$b_1\ne b'_1$ and $b_2\ne b'_2$ is true. Let $j\in\{1,2\}$ be the index
for which $|B_j|\ne|B'_j|$ (then $|B_{3-j}|=|B'_{3-j}|$).

\smallskip

\case 4{$\min\{|B_j|,|B'_j|\}\le n/3$}

\noindent
Spoiler wins in at most $\min\{|B_j|,|B'_j|\}+3\le n/3+3$ moves with
at most 1 alternation similarly to Case~3.

\smallskip

\case 5{$\min\{|B_j|,|B'_j|\} > n/3$}

\noindent
This is the only place in the proof where we need to assume
that $|G'|=|G|$. We can do so because in Case 5 the graph $G$
must have isolated vertices or edges. It follows that
$G[\compl B]$ and $G'[\compl{B'}]$
are non-isomorphic graphs of different orders.
Denote the graph of the smaller order by $H$ and the other
graph by $H'$.
Spoiler enforces play on $H$ and $H'$ as follows:
As soon as Duplicator moves in $B$ or $B'$, Spoiler, who has selected
in this round a vertex $u$ in a component $C$ of $G$ or $G'$ with at
least 3 vertices, wins in 2 extra moves by selecting two more vertices
$u_1$ and $u_2$ in $C$ so that $u$, $u_1$, and $u_2$ span a connected
subgraph. Note that
$$
|H|\le n-b<2n/3.
$$

As long as Duplicator moves outside $B$ and $B'$, Spoiler
uses the following strategy for the \EF\/ game on $H$ and $H'$.
We will now refer to Definition \ref{def:add}.
If $H'\ne H\add lv$ for any $v\in V(H)$, then Spoiler follows the
strategy from Lemma \ref{lem:d1} and wins in at most $|H|/2+2$
moves with at most 1 alternation.
If $H'=H\add lv$ for some $v\in V(H)$, then Spoiler follows the
strategy from Lemma \ref{lem:close} and wins in at most
$$
\frac{|H|+1}{\gs_H(v)+1}+\gs_H(v)+1\le\frac{|H|}2+\gs(H)+\frac32
$$
moves with at most 1 alternation between the graphs.
As $H$ has no isolated vertices, by Lemma \ref{lem:gsvsdelta}
we have $\gs(H)\le\Delta(H)+1\le d+1$.
Thus, in Case 5 Spoiler needs at most
$|H|/2+d+9/2< n/3+d+9/2$ moves to win.

\smallskip

In any of Cases 1--5 Spoiler makes at most 1 alternation and
at most $(\frac12-\frac12\tau)n+d^2+d+\frac72$ moves, which is actually
the claimed bound.

\medskip

\centerline{{\sc Strategy 2} (applicable if $a\le\tau n$)}

\smallskip

\noindent
We split our description of Spoiler's strategy into four phases.

\smallskip

{\it Phase 1.}

\noindent
Spoiler selects all vertices in the set $A$.

\smallskip

{\it Phase 2.}

\noindent
Spoiler will make moves in pairs. Let $i\ge1$.
Denote the vertices selected by him in the $(2i-1)$-th
and $2i$-th rounds of Phase 2 by $x_i$ and $y_i$ respectively.
Suppose that Spoiler has already made $2(i-1)$ moves and
selected a set $X_{i-1}=A\cup\{x_1,y_1,\dots,x_{i-1},y_{i-1}\}\subset V(G)$.
Let us explain how $x_i$ and $y_i$ are now selected.

If there is a vertex $x\in V(G)$ such that
\begin{itemize}
\item
$\dist(x,X_{i-1})\ge5$ and
\item
for any $y$ with $\dist(x,y)\le 2$ we have $\deg(y)\le \deg(x)$,
\end{itemize}
then Spoiler selects this $x$ for $x_i$.

\begin{claim}\label{cl:uyv}
Suppose that $x_i=x$ does exist. Then there are vertices $u,y,v$ such that
$\{x,u\},\{u,y\},\{y,v\}\in E(G)$ while $\{x,y\},\{x,v\}\notin E(G)$.
\end{claim}

\begin{subproof}
Let $C$ be the component of $G$ containing $x$.
It should contain a vertex $v$ with $\dist(x,v)=3$ for else
every vertex of $C$ would be at distance at most 2 from $x$
and hence $C$ would have at most $1+d+d(d-1)=d^2+1$ vertices.
Let $(x,u,y,v)$ be an arbitrary path from $x$ to $v$.
The vertices $u,y,v$ are as desired.
\end{subproof}

If $x_i=x$ is selected, Spoiler chooses some $u,y,v$ as in the claim
and takes the $y$ for~$y_i$.

If no such $x$ exists, Phase 2 ends.
Suppose that Phase 2 lasts $2r$ rounds.
Recall that the partition $\calC(X)$, where $X\subset V(G)$,
is defined by Definition \ref{def:calc}.

\begin{claim}\label{cl:cxi}
$|\calC(X_i)|\ge|\calC(X_{i-1})|+3$ if $i<r$ and
$|\calC(X_r)|\ge|\calC(X_{r-1})|+2$.
\end{claim}

\begin{subproof}
We will show that, if we extend $X_{i-1}$ to $X_i$, one of the classes
in $\calC(X_{i-1})$ splits up into at least 4 parts if $i<r$ and
into at least 3 parts if $i=r$.

By the choice of $u$, $y=y_i$, and $v$, we have $\dist(x_i,u)=1$
and $\dist(x_i,v)=3$.
Since $\dist(x_i,X_{i-1})\ge5$, neither $u$ and $v$ is in $X_{i-1}$.
Note that $u$ is adjacent to both $x_i$ and $y_i$, while
$v$ is adjacent to $y_i$ but not to $x_i$.

Note also that $\Gg(x_i)\setminus\Gg(y_i)\ne\emptyset$ because
$x_i$ has no less neighbors than $y_i$ has and $v$ is a neighbor
of $y_i$ but not of $x_i$. Thus, there is a vertex $w$ adjacent
to $x_i$ but not to $y_i$. Like $u$ and $v$, we have $w\notin X_{i-1}$.

Thus, $u$, $v$, and $w$ belong to pairwise distinct classes of $\calC(X_i)$.
If we assume that $i<r$, we are able to find a vertex in a yet another
class. Indeed, consider $x=x_{i+1}$. Since $\dist(x,X_i)\ge5$,
this vertex is adjacent neither to $x_i$ nor to $y_i$.

On the other hand, every of $u$, $v$, $w$, and $x$ is at distance
at least 2 from $X_{i-1}$. Therefore all of them are in the same
class of $\calC(X_{i-1})$.
\end{subproof}

{\it Phase 3.}

\noindent
As long as possible, Spoiler extends $X=X_r$ by one vertex so that
$|\calC(X)|$ increases at least by 1. Phase 3 ends as soon as Spoiler
arrives at a $\calC$-maximal set in the sense of Definition \ref{def:cmax}.

Suppose that Phase 3 lasts $h$ rounds. At the end of this phase
we therefore have
$$
|X|=a+2r+h
$$
and
$$
|\calC(X)|\ge 1+3(r-1)+2+h=3r+h.
$$
It follows that $|\calC(X)|\ge|X|+r-a$ and hence
$n\ge|X|+|\calC(X)|\ge2|X|+r-a$. We conclude that
$$
|X|\le\frac{n+a-r}2.
$$

{\it Phase 4.}

\noindent
Spoiler now plays precisely as in Phase 2 of the strategy
designed in the proof of Lemma \ref{lem:d1}.
By Lemma \ref{lem:usefulgs}, with Lemma \ref{lem:gsvsdelta}
taken into account, Spoiler wins making totally at most
\begin{equation}\label{eq:fin}
|X|+d+3\le\frac{n+a-r}2+d+3\le\frac{n+\tau n-r}2+d+3
\end{equation}
moves.
It therefore remains to show that the duration of Phase 2, controlled
by $r$, is linearly related to $n$ (the parameter $\tau$ is chosen
small enough).

\begin{claim}\label{cl:vk}
Let $V_k=\setdef{x\in V(G)}{\dist(x,X_r)\ge 2k+3}$.
Then $V_{d+1}=\emptyset$.
\end{claim}

\begin{subproof}
Assume, to the contrary, that $V_{d+1}\ne\emptyset$. We will show that
then there is $x_{r+1}$ such that $\dist(x_{r+1},X_r)\ge5$
and every $y$ at distance at most 2 from $x_{r+1}$ has smaller degree,
contradicting the fact that Phase 2 lasts $2r$ rounds.
Let $d_i=\max\setdef{\deg(x)}{x\in V_i}$.
If no $z_i\in V_i$ with $\deg(z_i)=d_i$ can be taken for $x_{r+1}$,
then $d_{i+1}<d_i$. Indeed, for any such $z_i$ there is $y$ such that
$\dist(z_i,y)\le2$ and $\deg(y)>\deg(z_i)$.
The latter implies that $y\notin V_i$, i.e., $\dist(y,X_r)\le 2i+2$.
It follows that $\dist(z_i,X_r)\le 2i+4$, hence $z_i\notin V_{i+1}$
and $d_{i+1}<d_i$. Since the chain $d_1>d_2>d_3>\ldots$ can have
length at most $d$, some $z_i$ with $i\le d+1$ can be taken for $x_{r+1}$.
This contradiction proves the claim.
\end{subproof}

Thus, $|\compl{V_{d+1}}|=n$. By the definition of $V_{d+1}$ we have
$$
|\compl{V_{d+1}}|\le|X_r|(1+d+d(d-1)+d(d-1)^2+\ldots+d(d-1)^{2d+3})
<|X_r|d^{2d+5}
$$
(note that $d\ge2$, which follows from the assumption that $a<n$).
Putting it together, under the assumption that $a\le \tau n$, we obtain
$$
\tau n+2r\ge a+2r=|X_r|>n/d^{2d+5},
$$
which implies that
$$
r>n(d^{-2d-5}-\tau)/2.
$$
Substituted in \refeq{eq:fin}, this shows that Spoiler wins in at most
$$
\of{\frac12+\frac34\tau-\frac14d^{-2d-5}}n+d+3
$$
moves, which is within the required bound.
\end{proofof}

\section{The worst case dimension of the \WL\/ algorithm}\label{s:wl}

The main purpose of this, mostly expository, section is to give
a self-contained proof of an upper bound for the dimension
of the \WL\/ algorithm for the graph isomorphism problem.
The dimension is an important parameter of the algorithm.
On the one hand, the higher dimension is chosen,
the longer the algorithm runs. On the other hand, a small dimension
may do not suffice to compute the right output.
Our goal will be to show that, on input graphs of order $n$,
the dimension $\lfloor(n+1)/2\rfloor$ suffices.
Another job we do here is to compute an explicit constant in the
Cai-F\"urer-Immerman bound~\refeq{eq:cfi}.

We begin with description of the algorithm.

\subsection{Definitions and notation}

Given an ordered $k$-tuple of vertices $\baru=(u_1,\ldots,u_k)\in V(G)^k$,
let $s=s(\baru)$ be the number of distinct components in $\baru$
and define a function $F_\baru\function{\{1,\ldots,k\}}{\{1,\ldots,s\}}$
by $F_\baru(i)=|\{u_1,\ldots,u_i\}|$.
Furhtermore, let $G_\baru$ be the graph on the vertex set $\{1,\ldots,s\}$
with vertices $a$ and $b$ adjacent iff, for the smallest $i$ and $j$
such that $F_\baru(i)=a$ and $F_\baru(j)=b$, $u_i$ and $u_j$ are
adjacent in $G$. The pair $(F_\baru,G_\baru)$ is an {\em isomorphism type\/}
of $\baru$ and will be denoted by $\isotype\baru$.

If $w\in V(G)$ and $i\le k$, we let $\baru^{i,w}$ denote the result
of substituting $w$ in place of $u_i$ in $\baru$.

A {\em (vertex) coloring\/} of a graph $G$ is an arbitrary function
$\gamma\function{V(G)}C$. We will say that $C$ is the {\em set of colors\/}
and that a vertex $v\in V(G)$ has color $\gamma(v)$. For a color $c\in C$,
the set $\gamma^{-1}(c)$ is its {\em monochromatic class}.
A coloring $\gamma'$ {\em refines\/} a coloring $\gamma$ if
$\gamma'(v)=\gamma'(u)$ implies $\gamma(v)=\gamma(u)$, that is,
the partition of $V(G)$ into $\gamma'$-monochromatic classes
refines the partition into $\gamma$-monochromatic classes.

\subsection{Description of the algorithm}

We distinguish two modes of the algorithm. In the {\em canonization
mode\/} the algorithm takes as an input a graph $G$ and is purported
to output its canonic form $W(G)$. The {\em canonic form\/} of a graph
$W$ is a graph function such that $W(G)=W(G')$ iff $G\cong G'$.
In the {\em isomorphism testing mode\/}
the algorithm takes as an input two graphs $G$ and $G'$ and should decide
if $G\cong G'$.

We now describe the $k$-dimensional algorithm. The algorithm
assigns an {\em initial coloring\/} to an input graph,
then step by step refines it by iterating the {\em color refinement\/}
procedure, and finally, when no color refinement is any more possible,
terminates and computes an output.

\medskip

{\sc Initial coloring}

\smallskip

\noindent
The algorithm assigns each $\baru\in V(G)^k$ color
$\wl Gk0\baru=\isotype\baru$ (in a suitable encoding).

\medskip

{\sc Color refinement step}

\smallskip

\noindent
In the $r$-th step each $\baru\in V(G)^k$ is assigned color
$$
\wl Gkr\baru=\of{\wl Gk{r-1}\baru,
\setdef{(\wl Gk{r-1}{\baru^{1,w}},\ldots,\wl Gk{r-1}{\baru^{k,w}})}%
{w\in V(G)}}.
$$
In the proper \WL\/ algorithm the second component of $\wl Gkr\baru$
is a multiset rather than a set. However, in what follows we assume
that it is a set. It will be clear that the version of the algorithm we
consider is weaker than the standard one, i.e.,
whenever our $k$-dimensional version gives the right output, so
does the standard $k$-dimensional version (and it will be not hard
to show that sometimes for the standard version a considerably smaller
dimension is enough). This relaxation makes our result only stronger
as any upper bound for the dimension of the weaker version is as well an
upper bound for the dimension of the standard version.

\begin{proposition}\label{prop:iso}
If $\phi$ is an isomorphism from $G$ to $G'$, then for all $k$, $r$,
and $\baru\in V(G)^k$ we have
$\wl Gkr\baru=\wl{G'}kr{\phi^k(\baru)}$. \noproof
\end{proposition}

\begin{proposition}\label{prop:stable}
For every pair of graphs $G$ and $G'$ there is a number $R$ such that
for all $\baru\in V(G)^k$, $\barv\in V(G')^k$, and $r>R$
$$
\wl Gkr\baru=\wl {G'}kr\barv\mbox{\ \ iff\ \ }
\wl GkR\baru=\wl {G'}kR\barv.
$$
Moreover, if $R_k(G,G')$ denotes the smallest such $R$, then
$R_k(G,G')<|G|^k+|G'|^k$.
\end{proposition}

\begin{proof}
By Proposition \ref{prop:iso} it suffices to prove the claim for
arbitrary isomorphic copies of $G$ and $G'$ and we therefore can
suppose that $V(G)$ and $V(G')$ are disjoint. Colorings $\wlcol Gkr$
and $\wlcol {G'}kr$ determine the partition of the union $V(G)^k\cup V(G')^k$
into monochromatic classes. Denote this partition by $\Pi^r$.
Since the $(r+1)$-th color incorporates the $r$-th color, $\Pi^{r+1}$
is a subpartition of $\Pi^r$. It is clear that we eventally have
$\Pi^{R+1}=\Pi^R$ and the smallest such $R$ is less than
$|V(G)|^k+|V(G')|^k$.
\end{proof}

\medskip

{\sc Computing an output}

\smallskip

\noindent
{\it Isomorphism testing mode.}
The algorithm terminates color refinement as soon as the partition
$\Pi^r$ of $V(G)^k\cup V(G')^k$ coincides with $\Pi^{r-1}$, i.e., after
performing $r=R_k(G,G')+1$ refinement steps. The algorithm decides that
$G\cong G'$ iff
\begin{equation}\label{eq:decision}
\setdef{\wl Gkr{u^k}}{u\in V(G)}=
\setdef{\wl{G'}kr{v^k}}{v\in V(G')},
\end{equation}
where $w^k$ denotes the diagonal vector $(w_1,\ldots,w_k)$ with all $w_i=w$.

{\it Canonization mode.}
The algorithm performs $r=2|G|^k-1$ refinement steps and outputs the set
$\setdef{\wl Gkr{u^k}}{u\in V(G)}$.

\paragraph{Implementation details and complexity bounds.}
Denote the minimum length of the code of $\wl Gkr\baru$ over all $\baru$
by $\lambda(r)$. As easily seen, for any natural encoding we should expect
that $\lambda(r)\ge(k+1)\lambda(r-1)$. To prevent increasing $\lambda(r)$
at the exponential rate,
before every refinement step we arrange colors of all $k$-tuples in the
lexicographic order and replace each color with its number.
In the canonization mode we should keep the substitution tables of all steps.
In the isomorphism testing mode this is unnecessary but it should be
stressed that color renaming must be common for both input graphs.
The straightforward implementation of the algorithm takes time
$O(k^2n^{2k}\log^2n)$ and space $O(kn^{2k}\log n)$, where $n=|G|$.
In the isomorphism testing mode, when we do not waste memory by
keeping substitution tables, the space $O(kn^k(k+\log n))$ suffices.
A better implementation with time bound $O(k^2n^{k+1}\log n)$ is
suggested in~\cite{ILa}.

\subsection{Relation to the \EF\/ game}

Given numbers $r$, $l$, and $k\le l$, graphs $G$, $G'$, and
$k$-tuples $\baru\in V(G)^k$,
$\barv\in V(G')^k$, we use notation $\game^l_r(G,\baru,G',\barv)$
to denote the $r$-round $l$-pebble \EF\/ game on $G$ and $G'$
with initial configuration $(\baru,\barv)$, i.e., the game starts
with one copy of the pebble $p_i$, $i\le k$, placed on $u_i$
and the other copy of $p_i$ placed on~$v_i$.

\begin{proposition}\label{prop:game}
{\bf (Cai, F\"urer, and Immerman \cite{CFI})}
For all $\baru\in V(G)^k$ and $\barv\in V(G')^k$ the equality
\begin{equation}\label{eq:coleq}
\wl Gkr\baru=\wl {G'}kr\barv
\end{equation}
holds iff Duplicator has a winning
strategy in $\game_r^{k+1}(G,\baru,G',\barv)$.
\end{proposition}

\begin{proof}
We proceed by induction on $r$. The base case $r=0$ is straightforward
by the definitions of the initial coloring and the game. Assume that
the proposition is true for $r-1$ rounds.

Assume \refeq{eq:coleq} and consider the \EF\/
game $\game_r^{k+1}(G,\baru,G',\barv)$.
First of all, the initial configuration is non-losing
for Duplicator since \refeq{eq:coleq} implies that
$\isotype\baru=\isotype\barv$. Further, Duplicator can survive in
the first round. Indeed, assume that Spoiler in this round selects
a vertex $a$ in one of the graphs, say in $G$. Then Duplicator selects
a vertex $b$ in the other graph, respectively in $G'$, such that
$\wl Gk{r-1}{\baru^{i,a}}=\wl {G'}k{r-1}{\barv^{i,b}}$ for all $i\le k$.
In particular,
$\isotype{\baru^{i,a}}=\isotype{\barv^{i,b}}$ for all $i\le k$.
Along with $\isotype\baru=\isotype\barv$, this implies that
$\isotype{\baru,a}=\isotype{\barv,b}$.
Assume now that in the second round Spoiler removes the $j$-th pebble,
$j\le k$. Then Duplicator's task in the rest of the game is
essentially to win $\game^{k+1}_{r-1}(G,\baru^{j,a},G',\barv^{j,b})$.
Since $\wl Gk{r-1}{\baru^{j,a}}=\wl {G'}k{r-1}{\barv^{j,a}}$,
Duplicator succeeds by the induction assumption.

Assume now that \refeq{eq:coleq} is false. It follows that
$\wl Gk{r-1}\baru\ne\wl{G'}k{r-1}\barv$ (then Spoiler has a winning
strategy by the induction assumption) or there is a vertex $a$ in one
of the graphs, say in $G$, such that for every $b$ in the other graph,
respectively in $G'$,
$\wl Gk{r-1}{\baru^{j_b,a}}\ne\wl {G'}k{r-1}{\baru^{j_b,b}}$
for some $j_b\le k$. In the latter case Spoiler in his first move
places the $(k+1)$-th pebble on $a$. Let $b$ be the vertex selected
in response by Duplicator. In the second move Spoiler will remove
the $j_b$-th pebble, which implies that since the second round
the players essentially play
$\game_{r-1}^{k+1}(G,\baru^{j_b,a},G',\barv^{j_b,b})$.
By the induction assumption, Spoiler wins.
\end{proof}

If $G\cong G'$, then the \WL\/ algorithm recognizes $G$ and $G'$ as
isomorphic for every dimension $k$. This follows from Proposition
\ref{prop:iso}. If $G\not\cong G'$, then the algorithm may be wrong
if $k$ is chosen too small.

\begin{corollary}\label{cor:wl}
If $G$ and $G'$ are non-isomorphic graphs of the same order $n$,
then the $k$-dimensional algorithm recognizes $G$ and $G'$ as
non-isomorphic iff $k\ge \V{G,G'}-1$.
\end{corollary}

\begin{proof}
Duplicator has a winning strategy in $\game_r^{k+1}(G,G')$ iff
for every $a\in V(G)$ (resp.\ $b\in V(G')$) there is
$b\in V(G')$ (resp.\ $a\in V(G)$) such that Duplicator has a winning
strategy in $\game_{r-1}^{k+1}(G,a,G',b)$ or, equivalently, in
$\game_{r-1}^{k+1}(G,a^k,G',b^k)$. It follows by Propositions
\ref{prop:stable} and \ref{prop:game} that the $k$-dimensional
algorithm decides that $G\cong G'$ iff Duplicator has a winning
strategy in $\game_r^{k+1}(G,G')$ for all $r$.
By Proposition \ref{prop:loggames},
$\V{G,G'}-1$ is equal to the maximum $l$ such that Duplicator
has a winning strategy in $\game_r^l(G,G')$ for every $r$.
Therefore the decision of the $k$-dimensional algorithm is correct
iff $k\ge \V{G,G'}-1$.
\end{proof}

\subsection{An upper bound on the dimension of the algorithm}

\begin{definition}
The smallest dimension of the \WL\/ algorithm giving the right output
on graphs $G$ and $G'$ of order $n$ in the isomorphism testing mode
will be referred to as the {\em optimum dimension\/} of the algorithm
on $G$ and $G'$ and denoted by $\Wl(G,G')$. Furthermore, the optimum
dimension of the algorithm on graphs of order $n$ is defined by
$$
\Wl(n)=\max\setdef{\Wl(G,G')}{G'\not\cong G,\ |G'|=|G|=n}.
$$
\end{definition}
It is easy to see that the smallest dimension of the algorithm giving
the right output on an input graph $G$ in the canonization mode
is equal to
$$
\max\setdef{\Wl(G,G')}{G'\not\cong G,\ |G'|=|G|}
$$
and hence in the worst case is equal to $\Wl(n)$.

On the account of Corollary \ref{cor:wl} and Theorem \ref{thm:d1},
we immediately obtain the following result.

\begin{theorem}\label{thm:wl}
$\Wl(n)\le(n+1)/2$.
\end{theorem}

This bound is almost tight for the relaxed version of the algorithm
that we have dealt with in this section.
However, Theorem \ref{thm:wl} leaves a considerable gap if compared
with the Cai-F\"urer-Immerman lower bound for $\Wl(n)$, that is
discussed in the next subsection.
Thus, what $\limsup_{n\to\infty}\Wl(n)/n$ is remains open.

\subsection{The lower bound for the dimension ---
computing an explicit constant}\label{ss:cfi}

Cai, F\"urer, and Immerman \cite{CFI} prove a striking linear lower bound
$\Wl(n)\ge cn$, without specification of a positive constant $c$.
We are curious to draw from their proof an explicit value of $c$.
We do not give a complete overview of the proof focusing only on
a few most relevant points. The following notion will be fairly useful.

\begin{definition}
Let $H$ be a graph of order $n$. Given $X\subset V(H)$, denote
$H\setminus X=H[V(H)\setminus X]$. We call a set $X$
a {\em separator\/} of $H$ if every connected component of the graph
$H\setminus X$ has at most $n/2$ vertices.
The number of vertives in a separator is called its {\em size}.
The minimum size of a separator of $H$ is denoted by $s(H)$.
\end{definition}

Cai, F\"urer, and Immerman present a construction of non-isomorphic
graphs $G$ and $G'$ of the same order with large $\Wl(G,G')$.
Both $G$ and $G'$ are constructed from a suitable connected graph $H$
with minimum vertex degree at least 2. We will assume that $H$ is
$d$-regular, that is, every its vertex has degree $d$
(using $H$ not regular seems to give us no gain).
Below we summarize the properties of the construction.
Recall that $\Delta(G)$ denotes the maximum vertex degree of a graph~$G$.

\begin{proposition}\label{prop:cfi}
{\bf (Cai, F\"urer, and Immerman \cite{CFI})}
Let $d\ge2$ and $H$ be a connected $d$-redular graph.
There are transformations of $H$ to two graphs
$G=G(H)$ and $G'=G'(H)$ such that
\begin{itemize}
\item
If $d\ge3$, both $G$ and $G'$ are connected;
\item
$|G|=|G'|=(d+2^{d-1})|H|$;
\item
$\Delta(G)=\Delta(G')=2^{d-1}$;
\item
$G\not\cong G'$;
\item
$\Wl(G,G')\ge s(H)$.\noproof
\end{itemize}
\end{proposition}
Thus, we need a family of $d$-regular graphs $H$ with $d$ constant
and $s(H)$ linearly related to the order of $H$. The authors of $\cite{CFI}$
suggest using graphs with good expansion properties.

\begin{definition}
Let $H$ be a graph of order $n$. The {\em vertex-expansion\/}
of $H$ is denoted by $i_v(H)$ and defined by
$$
i_v(H)=\min\setdef{\frac{|N(A)|}{|A|}}{A\subset V(H),\ |A|\le\frac n2},
$$
where $N(A)=\bigcup_{v\in A}\Gg(v)\setminus A$ is the {\em neiborhood\/}
of a set~$A$.
\end{definition}

\begin{lemma}\label{lem:expan}
For a graph $H$ of order $n$ we have
$$
s(H)\ge\frac{i_v(H)}{3+i_v(H)}\,n.
$$
\end{lemma}

\begin{proof}
Let $X$ be a separator of $H$ with the smallest size $s=s(H)$.
Denote the largest size of a connected component of $H\setminus X$
by $m$ and recall that $m\le n/2$.
There is a set $A_1\subseteq V(G)\setminus X$ with $|A_1|=m$ such
that $H[A_1]$ is a connected component of $H\setminus X$ and, as
it is not hard to see, there is a set $A_2\subseteq V(G)\setminus X$
with $\frac{n-s}2-\frac m2\le|A_2|\le\frac{n-s}2$ such
that $H[A_2]$ is a union of connected components of $H\setminus X$.
Note that $\max\{m,\frac{n-s}2-\frac m2\}\ge\frac{n-s}3$ and
therefore for $A_i$, one of the sets $A_1$ and $A_2$, we have
$$
\frac{n-s}3\le|A_i|\le\frac n2.
$$
By the definition of the vertex expansion,
$$
|N(A_i)|\ge i_v(H)|A_i|\ge i_v(H)\,\frac{n-s}3.
$$
Since $A_i$ is a union of connected components of $H\setminus X$,
we have $N(A_i)\subseteq X$ and hence $|N(A_i)|\le s$.
Thus, we obtain the relation
$$
s\ge i_v(H)\,\frac{n-s}3.
$$
Resolving it with respect to $s$, we arrive at the reqiured estimate.
\end{proof}

Thus, we need $d$-regular graphs with large vertex-expansion.
The best examples we could find in the literature come from the
known edge-expansion results.

\begin{definition}
Let $H$ be a graph of order $n$. The {\em edge-expansion\/}
(or the {\em isoperimetric number\/}) of $H$ is denoted by
$i_e(H)$ and defined by
$$
i_e(H)=\min\setdef{\frac{e(A,N(A))}{|A|}}{A\subset V(H),\ |A|\le\frac n2},
$$
where $e(A,B)$ denotes the number of edges in $H$ with one
end vertex in $A$ and the other in~$B$.
\end{definition}

For a $d$-regular graph $H$ it is straightforward that $i_v(H)\ge i_e(H)/d$.
We are able to improve this relation for 3-regular (or {\em cubic\/}) graphs.

\begin{lemma}\label{lem:cubic}
If $H$ is a cubic graph, then $i_v(H)\ge i_e(H)/2$.
\end{lemma}

\begin{proof}
If $H$ is disconnected, then one of its connected components
occupies no more than a half of the vertices and hence $i_v(H)=0$.
Suppose that $H$ is connected.

Of all $A\subset V(H)$ with $|A|\le|H|/2$ and
$|N(A)|/|A|=i_v(H)$, take one minimizing $e(A,N(A))$.
Let us show that every vertex $x\in N(A)$ sends at most 2 edges to
$A$. This will give us the desired relation because in this case
$$
i_e(H)\le\frac{e(A,N(A))}{|A|}\le\frac{2|N(A)|}{|A|}=2\,i_v(H).
$$
Suppose, to the contrary, that some $x\in N(A)$ sends 3 edges to~$A$.

Consider an arbitrary $y\in A$. Let $A_y=(A\setminus\{y\})\cup\{x\}$.
If $y$ sends 3 edges to $N(A)\setminus\{x\}$, then
$N(A_y)=N(A)\setminus\{x\}$ has less vertices than $N(A)$ has
while $|A_y|=|A|$. Therefore $|N(A_y)|/|A_y|<i_v(H)$, a contradiction.
If $y$ sends 1 or 2 edges to $N(A)\setminus\{x\}$, then
$N(A_y)\subseteq (N(A)\setminus\{x\})\cup\{y\}$. It follows that
$|N(A_y)|\le|N(A)|$ and it should be $|N(A_y)|/|A_y|=i_v(H)$.
However, $e(A_y,N(A_y))\le e(A,N(A))-1$, contradicting the choice of~$A$.

We conclude that every $y\in A$ sends no edges to $N(A)\setminus\{x\}$.
It follows that $A\cup\{x\}$ spans a proper connected component of $H$,
a contradiction to the connectivity of~$H$.
\end{proof}

Using Lemmas \ref{lem:expan} and \ref{lem:cubic}, from
Proposition \ref{prop:cfi} we easily obtain the following consequence.

\begin{proposition}\label{prop:expl}
Let $i_e(3,m)$ denote the maximum edge-expansion of a connected cubic
graph of order $m$. Then there are non-isomorphic graphs
$G$ and $G'$ both of order $7m$ with maximum degree 4 such that
$$
\Wl(G,G')\ge \frac{i_e(3,m)}{6+i_e(3,m)}\,m.\ \bull
$$
\end{proposition}

It seems that the best known lower bounds for $i_e(3,m)$ are
obtained by examining random cubic graphs. The edge-expansion
of a random cubic graph was studied by Buser \cite{Bus},
Bollob\'as \cite{Bol}, and others
with the best lower bound as follows.

\begin{proposition}
{\bf (Kostochka and Melnikov \cite{KMe1,KMe2})}
Let $H$ be a random cubic graph of order $m$.
If $m$ is sufficiently large, then with probability $1-o(1)$
we have
$$
i_e(H)\ge \frac{1}{4.95}.\ \bull
$$
\end{proposition}

It follows that $i_e(3,m)\ge\frac{1}{4.95}$, where $m$ is supposed
to be large enough.
For the graphs $G$ and $G'$ as in Proposition \ref{prop:expl}
we therefore have
$$
\Wl(G,G')\ge\frac n{7(1+6\cdot 4.95)}> 0.00465\,n,
$$
where $n=7m$ is the order of the graphs. Thus, the constant in question
is evaluated.

\begin{proposition}
$\Wl(n)> 0.00465\,n$ for infinitely many $n$.\noproof
\end{proposition}

Notice that, with high probability, a random cubic graph is connected.
The construction of \cite{CFI} together with the logical characterization
of $\Wl(G,G')$ given in the same paper, has therefore the following
consequence worthwhile to note.

\begin{proposition}\label{prop:connbdeg}
For infinitely many $n$, there are non-isomorphic connected graphs
$G$ and $G'$ both of order $n$ with maximum degree 4 such that
$\D{G,G'}> 0.00465\,n$.\noproof
\end{proposition}

\section{Digraphs and binary structures}\label{s:digraphs}

We now extend the results of Sections \ref{s:dist} and \ref{s:def}
over directed graphs and, more generally, structures with relation
arity at most~2.

\subsection{Definitions}\label{ss:didef}

In logical terms, a {\em directed graph\/} (or {\em digraph\/}) $G$
is an arbitrary binary relation $E$ on a vertex set $V(G)$.
The {\em edge set\/} of $G$ is $E(G)=\setdef{(u,v)\in V(G)^2}{E(a,b)=1}$.
Thus, between two distinct vertices $u$ and $v$
we allow two opposite edges $(u,v)$ and $(v,u)$, or only one of them,
or none. We view an edge $(u,v)$ as an arrow from $u$ to $v$.
An edge $(v,v)$, called a {\em loop}, is also allowed.
{}From now on the stand-alone term {\em graph\/} means ordinary
undirected graph. Recall that the undirected graphs are actually
considered a subclass of the directed graphs
wherein an undirected edge $\{u,v\}$ corresponds to two directed edges
$(u,v)$ and $(v,u)$ and we have no loops and no two vertices with
exactly one directed edge between them.

A {\em (vertex) colored digraph\/} is a structure that, in addition to the
binary relation, has unary relations $U_1,\ldots,U_m$. The truth of $U_i(v)$
for a vertex $v$ is interpreted as coloration of $v$ in color $i$.
Thus, a vertex can have several colors or no color.
However, colored digraphs can be modelled as digraphs with
each vertex having exactly one color by defining new colors
as conjunctions of some of $U_1,\ldots,U_m$ and yet another new color
for uncolored vertices.

\smallskip

{\bf Convention.}
Observe that digraphs can be modelled as colored loopless digraphs by
assigning a special color to the vertices with loops.
To facilitate the exposition, we will use this observation and
assume throughout this section all digraphs loopless.

\smallskip

A {\em complete digraph\/} has two directed edges $(u,v)$ and
$(v,u)$ for every pair of distinct vertices $u$ and $v$.
An {\em empty digraph\/} has no edges at all.
Let $G$ be a digraph and $X,Y\subseteq V(G)$ be disjoint.
Similarly to graphs, $G[X]$ denotes the subdigraph induced by $G$ on $X$
and $G[X,Y]$ denotes the bipartite subdigraph induced on the
vertex classes $X$ and $Y$.
A set $X$ is called {\em complete\/} (resp.\
{\em independent}) if $G[X]$ is complete (resp.\ empty).
$X$ is called {\em homogeneous\/} if it is complete or empty.
A bipartite subdigraph $G[X,Y]$ is {\em complete}, {\em independent}, or
{\em $(X,Y)$-complete} if for any $u\in X$ and $v\in Y$ we have
respectively $(u,v),(v,u)\in E(G)$; $(u,v),(v,u)\not\in E(G)$;
$(u,v)\in E(G)$ but $(v,u)\not\in E(G)$. It is
{\em dicomplete\/} if it is either $(X,Y)$-complete or $(Y,X)$-complete.
A pair $X,Y$ is called {\em homogeneous\/} if $G[X,Y]$ is complete,
independent, or dicomplete.

A {\em binary vocabulary\/} has relation arities at most 2.
A {\em binary structure\/} is a structure over a binary vocabulary.
Combinatorially, a binary structure
$(U_1,\ldots,U_s,E_1,\allowbreak\ldots,E_t)$
with $s$ unary and $t$ binary relations can be viewed as a complete
digraph that is vertex colored in colors $1,\ldots,s$ and edge colored
in colors $1,\ldots,t$: A vertex $v$ has color $i\le s$ iff $U_i(v)=1$
and an edge $(u,v)$ has color $j\le t$ iff $E_j(u,v)=1$.
An isomorphism between such complete digraphs should preserve
the sets of colors of each vertex and of each edge.
Given a binary structure $(U_1,\ldots,U_s,E_1,\ldots,E_t)$,
we will often call elements of its universe vertices and
consider each $E_i$ a digraph.

Let $G=(U_1,\ldots,U_s,E_1,\ldots,E_t)$ be a binary structure
with universe $V(G)$. A set $X\subseteq V(G)$ is {\em monochromatic\/} if
$U_i(u)=U_i(v)$ for all $u,v\in X$ and every $i\le s$.
We call $X$ {\em homogeneous\/} if it is monochromatic and
homogeneous in every digraph $E_j$, $j\le t$.
We call a pair $X,Y\subseteq V(G)$ of disjoint sets {\em homogeneous\/}
if it is such in every digraph $E_j$, $j\le t$.

\subsection{Distinguishing binary structures}

We now generalize Theorem \ref{thm:d1} from graphs to binary structures.
The same proof in essence goes through with only a few substantiate
modifications; We indicate these by tracing Section \ref{s:dist}.
By $G$ we will mean, unless stated otherwise, a binary structure
$(U_1,\ldots,U_s,E_1,\ldots,E_t)$.

Definition \ref{def:sim}, as well as the other definitions in
Section \ref{ss:sprel}, makes a perfect sense for an arbitrary
structure. Let us look what the similarity of vertices $u$ and $v$
means in a digraph. We say that a vertex $t$ {\em separates\/} $u$ and $v$
if precisely one of $(t,u)$ and $(t,v)$ is in $E(G)$ or precisely one
of $(u,t)$ and $(v,t)$ is in $E(G)$. Then $u$ and $v$ are similar if
no third vertex separates them and if $(u,v)$ and $(v,u)$ both are
in $E(G)$ or both are not\footnote{and
if $u$ and $v$ simultaneously make loops or do not
(see Convention in Section \ref{ss:didef}).}.
In a general binary structure, $u$ and $v$ are similar
if they have the same sets of colors (i.e., satisfy the same
unary predicates) and are similar in each digraph $E_i$.

\begin{claim}
Lemma \ref{lem:simprop} holds true for digraphs and, more
generally, for binary structures.
\end{claim}

\begin{subproof}
The only not completely obvious part of the proof
is verification that $\sim$ is a transitive relation on the vertex set
of a digraph. Suppose that $u\sim v$ and $v\sim w$ for pairwise distinct
$u$, $v$, and $w$.
Assume, for example, that $(u,v)$ and $(v,u)$ are both present in $E(G)$.
As $v\sim w$, we have $(u,w),(w,u)\in E(G)$. Now, as $u\sim v$, we
have $(v,w),(w,v)\in E(G)$. Also the sets of in- and out-neighbors
of $u$ and $w$ in
$V(G)\setminus\{u,v,w\}$ are identical (being equal to those of $v$).
This implies that $u\sim w$, proving the transitivity.
\end{subproof}

\begin{claim}\label{cl:lemalpha}
Lemma \ref{lem:alpha} holds true for digraphs and, more
generally, for binary structures with one stipulation in Item 1.
Specifically, let $X\subset V(G)$ be $\calC$-maximal.
Then the partition $\calC(X)$ has the following properties.
\begin{enumerate}
\item
Every $C$ in $\calC(X)$ is a homogeneous set provided $|C|\ne2$.
\item
If $C_1$ and $C_2$ are distinct classes in $\calC(X)$ and have at least
two elements each, then the pair $C_1,C_2$ is homogeneous.
\end{enumerate}
\end{claim}

\begin{subproof}
The claim easily reduces to its particular case that $G$ is a digraph.
So we assume the latter.

1) Suppose, to the contrary, that $C$ is not homogeneous.

First, suppose there are $u,v\in C$ with $(u,v)\in E(G)$ and
$(v,u)\not\in E(G)$. Let $w\in C\setminus\{u,v\}$. As moving $v$ to
$X$ cannot separate $u$ and $w$ (for otherwise $|\calC(X)|$
increases), we conclude that $(w,v)\in E(G)$ and
$(v,w)\not\in E(G)$. Pondering $w$ as the next move we conclude that
$(w,u)\in E(G)$ and $(u,w)\not\in E(G)$. But $u$ separates $v$ and
$w$, a contradiction.

Thus we can assume that $C$ contains three vertices $u$, $v$, and
$w$ such that $(u,v),(v,u)\in E(G)$ but $(v,w),(w,v)\not\in E(G)$.
However, if we move $v$ to $X$, then $C$ splits into two classes,
which are non-empty as they contain $u$ and $w$ respectively. Hence
$|\calC(X)|$ increases at least by 1, a contradiction.

2)
Let $u,v\in C_1$ and $w,x\in C_2$. The contemplation of $w$ for the
next move shows that $(u,w)$ and $(v,w)$ are both present or absent;
considering $v$ we conclude the same about $(v,x)$ and $(v,w)$. Hence,
$(u,w)$ is an edge if and only if $(v,x)$ is. Similarly, the same is
true about $(w,u)$ and $(x,v)$. As $u,v,w,x$ are arbitrary, the claim
follows.
\end{subproof}

Lemma \ref{lem:beta} carries over with literally the same proof.

Definition \ref{def:add} makes a perfect sense for binary structures.
We now state an analog of Lemma~\ref{lem:d1}.

\begin{lemma}\label{lem:did1}
Suppose that $G$ and $G'$ are non-isomorphic structures over the same
vocabulary with maximum relation arity 2.
If $G$ and $G'$ have orders $n$ and $n'$ respectively and $n\le n'$, then
\begin{equation}\label{eq:did1n4}
\DD1{G,G'}\le(n+5)/2
\end{equation}
unless $G'=G\add(n'-n)v$ for some $v\in V(G)$.
\end{lemma}

\begin{proof}
We go through the lines of the proof of Lemma \ref{lem:d1}.
Phase 1 is played with no changes. The same strategy for Phase 2
ensures that, unless Spoiler wins in 2 extra moves, we have
a partial isomorphism $\phi^*\function{X\cup Y}{X'\cup Y'}$
from $G$ to $G'$. Exactly as in the case of graphs, in Phase 2
Duplicator should obey $\phi^*$ if Spoiler moves in $Y\cup Y'$
and the $\phi^*$-similarity if Spoiler moves in $Z\cup Z'$.
In particular, Claim \ref{cl:phi*} carries over literally
and we again have the $\phi^*$-similarity of the classes $D_i$
and $D'_i$ for all $i\le p$.
Claim \ref{cl:didj} needs a careful revision, after which it is provable
with minor changes.

\begin{claim}\label{cl:dididj}
Unless Spoiler is able to win making 2 moves and at most 1 alternation
in Phase 2, the following conditions are met for every $i\le p$,
for every pair of distinct $i,j\le p$, and for each $k\le t$.
\begin{enumerate}
\item
Both $D_i$ and $D'_i$ are monochromatic and,
moreover, the sets of their colors coincide.
\item
If $|D_i|\ne2$, then the set $D_i$ is complete or independent in
the digraph $E_k$. Moreover, irrespectively of $|D_i|$,
the set $D_i$ is complete (resp.\ independent) in $E_k$ iff so is $D'_i$
in $E'_k$.
\item
The pairs $D_i,D_j$ and $D'_i,D'_j$
are simultaneously complete, independent, or dicomplete in
$E_k$ and $E'_k$ respectively. Moreover, in case of dicompleteness we
have the same direction of edges.\qed
\end{enumerate}
\end{claim}

Assume that the conditions in Claim \ref{cl:dididj} are obeyed.
We now claim that $G$ and $G'$ are isomorphic
iff $|D_i|=|D'_i|$ for all $i\le p$. We still cannot substantiate this,
as it was done in the case of graphs,
because of the stipulation that $|D_i|\ne2$ made in Item 2 of
Claim \ref{cl:dididj}.
The next claim gives the missing part of the argument.

\begin{claim}\label{cl:did2}
Unless Spoiler is able to win making 2 moves and at most 1 alternation
in Phase 2, the following is true for every $i\le p$. Assume that
$|D_i|=2$ and, for some $k\le t$, the digraph $E_k[D_i]$
has exactly one directed edge. Then $|D'_i|=2$ and $\phi^*$ extends to
an isomorphism from $G[X\cup Y\cup D_i]$ to $G'[X'\cup Y'\cup D'_i]$.
\end{claim}

\begin{subproof}
Let $C$ be the class in $\calC(X)$
including $D_i$ and $C'$ be the class in $\calC(X')$ including $D'_i$.
Claim \ref{cl:lemalpha} implies that $C=D_i$.

Let us prove that $|D'_i|=2$. Note first that $D'_i$ cannot consists
of a single element. Indeed, since $|C'|\ge2$, this would mean that
$C'$ is split up into at least two $D'$-classes. However, at most
one of them can have $\phi^*$-similar counterpart, contradicting
the assumption that the $\phi^*$-similarity is a perfect matching
between $\calD(X)$ and $\calD(X')$.

Let us now show that $D'_i$ cannot have more than 2 elements.
Suppose, to the contrary, that $|D'_i|>2$. Let $D_i=\{u,v\}$
with $(u,v)\in E_k$ and $(v,u)\notin E_k$.
If $E'_k[C']$ contains two vertices with
both or no directed edges between them, Spoiler selects these two and
wins because Duplicator is enforced
to reply with $u$ and $v$ by the analog of Claim \ref{cl:phi2}.
Assume that between any two vertices of $E'_k[C']$ there is exactly
one directed edge. It is easy to see that for some $a,b,c\in C'$ we have
$(a,b),(b,c)\in E'_k$. Let Spoiler select $b$ first.
By the analog of Claim \ref{cl:phi2}, Duplicator must
reply with $u$ or $v$. In the next move Spoiler wins with $a$
if Duplicator replies with $u$, and wins with $c$ otherwise.

To prove that $\phi^*$ extends as desired,
let Spoiler select the two vertices of $D'_i\subseteq C'$.
By the analog of Claim \ref{cl:phi2}, Duplicator must
reply with the vertices of $D_i=C$. It is clear that
either Spoiler wins or $\phi^*$ does extend.
\end{subproof}

Since $G$ and $G'$ are supposed non-isomorphic, there is $D_i$
such that $|D_i|\ne|D'_i|$. As in the case of graphs, we call such
a $D_i$ {\em useful}. Note that, by Item 2 of Claim \ref{cl:dididj}
and Claim \ref{cl:did2}, every useful $D_i$ is homogeneous
and, by Claim \ref{cl:dididj}, $D'_i$ is also homogeneous coherently
with $D_i$ over all unary and binary relations.

The rest of the proof carries over without any changes.
\end{proof}

\noindent
The main result of this section is provable
similarly to Theorem \ref{thm:d1} virtually with no change.

\begin{theorem}\label{thm:did1}
Let $G$ and $G'$ be structures over the same vocabulary with
maximum relation arity 2.
If $G$ and $G'$ are non-isomorphic and have the same order $n$, then
$\DD1{G,G'}\le(n+3)/2$.
\end{theorem}

\subsection{Defining binary structures}

Theorem \ref{thm:def1} generalizes to binary structures with minor
changes in the proof as follows.
Lemmas \ref{lem:sim} and \ref{lem:close} hold
true for binary structures literally, see Remark \ref{rem:arbstr}
where $k=2$ makes no change. As a consequence, by Lemma \ref{lem:did1},
Lemma \ref{lem:upper} literally holds as well. Lemma \ref{lem:exact},
where the homogeneousity is redefined in Section \ref{ss:didef},
also holds literally with the same in essence, easily adapted proof.
Note that Definition \ref{def:cls} makes a perfect sense for
the generalized homogeneousity. Finally, Theorem \ref{thm:def1}
holds true with literally the same proof.

\section{Uniform hypergraphs}\label{s:hyper}

We here extend the results of Sections \ref{s:dist} and \ref{s:def}
to uniform hypergraphs.

\subsection{Definitions}

A {\em $k$-uniform hypergraph\/} (or {\em $k$-graph\/}) $G$
on vertex set $V(G)$ is a family of $k$-element subsets of $V(G)$
called {\em (hyper)edges}. As usually, the set of edges of $G$
is denoted by $E(G)$. The $k$-graphs generalize the notion of
an ordinary graph, which is actually a $2$-graph.
{}From the logical point of view, a $k$-graph
is a structure with a single $k$-ary relation $E$ which is symmetric
in the sense that $E(x_1,\ldots,x_k)=E(x_{\pi(1)},\ldots,x_{\pi(k)})$
for any permutation $\pi$ of the $k$-element index set and
anti-reflexive in the sense that $E(x_1,\ldots,x_k)=0$ whenever
at least two of the $x_i$'s coincide. Thus, we will sometimes write
$E(x_1,\ldots,x_k)=1$ to say that $\{x_1,\ldots,x_k\}\in E(G)$
and $E(x_1,\ldots,x_k)=0$ to say that $\{x_1,\ldots,x_k\}\notin E(G)$.

Most graph-theoretic notions and all notions introduced for
general structures directly extend to $k$-graphs. For example,
the {\em complete\/} $k$-graph of order $n$ has $n$ vertices and
all possible ${n\choose k}$ edges. The {\em empty\/} $k$-graph of
order $n$ has $n$ vertices and no edges. If $X\subseteq V(G)$, then
$G[X]$ denotes the $k$-graph induced by $G$ on $X$, that is, the
$k$-graph with the vertices in $X$ and with all those edges $A\in E(G)$
for which $A\subseteq X$.

Let $0\le b\le k-1$ and $a=k-b$.
Given a $b$-vertex set $B\subset V(G)$, we define
the {\em link-graph\/} $G_B$ of $B$ to be the $a$-graph on
the vertex set $V(G)\setminus B$ with the edge set
$$
E(G_B)=\setdef{A}{|A|=a,\ A\cup B\in E(G)}.
$$
Clearly, $G_\emptyset=G$.

\subsection{Distinguishing $k$-graphs}

\begin{theorem}\label{thm:hd}
Let $k\ge2$.
If $G$ and $G'$ are non-isomorphic $k$-graphs both of order $n$, then
$$
\DD1{G,G'}\le\of{1-\frac1k}n+2k-1.
$$
\end{theorem}

If $k>2$, we do not know if the bound of the theorem is tight
since the only lower bound we know for any $k$ is $(n+1)/2$.
The latter is given by the same Example \ref{ex:lower} as for
$2$-graphs, where $K_m$ now means the complete $k$-graph and
$\compl{K_m}$ means the empty $k$-graph of order~$m$.

The proof of Theorem \ref{thm:hd} takes the rest of the subsection.
It will be built on the framework worked out in Section \ref{s:dist}.
However, note that Spoiler's strategy that will be designed in the
proof of the key Lemma \ref{lem:hd1} will be somewhat different.

As it was already mentioned, Definition \ref{def:sim} of the similarity
relation $\sim$ makes a perfect sense for an arbitrary class of
structures, in particular, for $k$-graphs.
Recall that $u\sim v$ for vertices $u$ and $v$ of a $k$-graph $G$
if the transposition $(uv)$ is an automorphism of $G$.
Both items of Lemma \ref{lem:simprop} hold true for $k$-graphs
but Item 2 should be supplemented with more information specific for
hypergraphs, see Item 4 of the following lemma.

\begin{lemma}\label{lem:hsimprop2}
Let $G$ be a $k$-graph.

\begin{enumerate}
\item
Given a $k$-element $U\subseteq V(G)$ and $v,w\in V(G)$,
let $U^{(vw)}$ denote the result
of substituting $v$ in place of $w$ and vise versa in $U$.
Assume that $v\sim w$. Then $U\in E(G)$ iff
$U^{(vw)}\in E(G)$.
\item
${\sim}$ is an equivalence relation on $V(G)$.
\item
If $v_i\sim w_i$ for all $1\le i\le k$, then
$E(v_1,\ldots,v_k)=E(w_1,\ldots,w_k)$.
\item
Let $v\in V(G)$ and $[v]_G$ denote the similarity class containing $v$.
Let $B\subset V(G)$ be disjoint with $[v]_G$ and have at most $k-1$ vertices.
Then the graph $G_B\left[[v]_G\right]$ is either complete or empty.
In particular, $G\left[[v]_G\right]$ is either complete or empty.
\item
Let $H$ be another $k$-graph, $X\subseteq V(H)$, and $G=H[X]$.
Then $[v]_H\cap X\subseteq[v]_G$ for any $v\in X$. In other words,
the partition of $X$ into the $H$-similarity classes refines the
partition into $G$-similarity classes.
\end{enumerate}
\end{lemma}

\begin{proof}
1)
This item is straightforward from the definition of the $\sim$-relation.

2)
The only not completely obvious task is, given pairwise distinct
vertices $u$, $v$, and $w$, to infer from $u\sim v$ and $v\sim w$
that $u\sim w$. For any $(k-1)$-element set of vertices
$\{x_1,\ldots,x_{k-1}\}$ we have to check that
\begin{equation}\label{eq:trans}
E(u,x_1,\ldots,x_{k-1})=E(w,x_1,\ldots,x_{k-1}).
\end{equation}
If $x_i\ne v$ for any $i$, this is easy. Assume that $x_{k-1}=v$.
Then the relations $u\sim v$ and $v\sim w$ imply by Item 1 that both
the left hand side and the right hand side of \refeq{eq:trans} are
equal to $E(u,w,x_1,\ldots,x_{k-2})$.

3)
Let $t$ denote the number of common elements in $\{v_1,\ldots,v_k\}$
and $\{w_1,\ldots,w_k\}$. We proceed by the reverse induction on $t$.
If $t=k$, the claim is trivial. Otherwise assume that
$w_i\notin\{v_1,\ldots,v_k\}$. We have
$E(v_1,\ldots,v_i,\ldots,v_k)=E(v_1,\ldots,w_i,\ldots,v_k)$
by Item 1 and
$E(v_1,\ldots,w_i,\ldots,v_k)=E(w_1,\ldots,w_i,\ldots,w_k)$
by the induction hypothesis.

4)
This item follows directly from Item 3.

5)
Let $v,w\in X$ and assume that $v\sim w$ in $H$.
As easily follows from Item 1, $v\sim w$ as well in~$G$.
\end{proof}

We will use all the definitions and notation introduced
in Section \ref{ss:sprel}. In particular, given $X\subset V(G)$,
we will deal with the relation $\xeq$ on $\compl X$, the partition
$\calC(X)$ of $\compl X$, the sets $Y(X)$ and $Z(X)$, and
the partition $\calD(X)$ of $Z(X)$. In addition, we need a new
relation and a new important parameter.

\begin{definition}
The relation ${\simeq}_X$ coincides with the equality on $X\cup Y(X)$
and with the relation $\xyxeq$ on $Z(X)$.
\end{definition}

\begin{definition}
$\tau(X)=|Y(X)|+|\calD(X)|$.
\end{definition}

We now generalize Definition \ref{def:add} over $k$-graphs.
This needs more care because the operation $\add$ for $k$-graphs
is somewhat more subtle than for graphs.

\begin{definition}
Let $G$ and $H$ be $k$-graphs, $v\in V(G)$, and $l\ge0$.
The notation $H\doteq G\add lv$ means that the following conditions
are fulfilled.
\begin{description}
\item[A1]
$\gs_G(v)\ge k$;
\item[A2]
$|H|=|G|+l$ and $V(G)\subseteq V(H)$;
\item[A3]
$H[V(G)]=G$;
\item[A4]
$[v]_H=[v]_G\cup(V(H)\setminus V(G))$.
\end{description}
Furthermore, we write $H=G\add lv$ if there is a $k$-graph $K$ such that
$H\cong K$ and $K\doteq G\add lv$.
\end{definition}

\noindent
Let us see carefully what the relation $H\doteq G\add lv$ means for
$k$-graphs. As easily seen, $H$ is obtained from $G$ by adding a set $A$
of $l$ new vertices and some new edges involving at least one vertex
from $A$. Given a $k$-vertex set $U\subseteq V(H)$ having non-empty
intersection with $A$, we have to decide whether or not $U$
is in $E(H)$. The criterion is given by Item 4 of Lemma \ref{lem:hsimprop2}.
Specifically, assume that the intersection $B=U\cap(V(G)\setminus[v]_G)$
contains $b$ vertices. Let $W\subseteq[v]_G$ with $|W|=k-b$.
Then $U\in E(H)$ iff $B\cup W\in E(G)$ (the latter is equally true
or false for any choice of $W$). This provides us with the complete
description of $H$. Thus, on the account of Item 4 of
Lemma \ref{lem:hsimprop2}, we arrive at the conclusion that
$H=G\add lv$, if exists, is unique up to an isomorphism.
It remains to prove that, if $H$ is constructed as described above,
then indeed $H\doteq G\add lv$. This immediately follows from the
following lemma.

\begin{lemma}\label{lem:addequi}
Let $G$ and $H$ be $k$-graphs, $v\in V(G)$, and $l\ge0$.

\smallskip

\noindent
{\rm 1)}\ \
$H\doteq G\add lv$ iff the following two conditions are fulfilled.
\begin{description}
\item[B1]
$|H|=|G|+l$ and $V(G)\subseteq V(H)$;
\item[B2]
There is $D\subseteq[v]_G$ with $|D|\ge k$ such that the following is true:
Every injection $\psi\function{V(G)}{V(H)}$ whose restriction to
$V(G)\setminus D$ is the identity map is a partial isomorphism from
$G$ to $H$.
\end{description}

\noindent
{\rm 2)}\ \
$H=G\add lv$ iff the following two conditions are fulfilled.
\begin{description}
\item[C1]
$|H|=|G|+l$;
\item[C2]
There is $D\subseteq[v]_G$ with $|D|\ge k$ such that the following is true:
There exists an injection $\psi_0\function{V(G)\setminus D}{V(H)}$
whose every injective extension $\psi\function{V(G)}{V(H)}$
is a partial isomorphism from $G$ to $H$.
\end{description}
\end{lemma}

\begin{proof}
1)
Let us show first that Conditions A1--A4 imply Conditions B1--B2.
Since B1 coincides with A2, we focus on B2.
Let $D\subseteq[v]_G$ be an arbitrary set with $|D|\ge k$,
existing by A1. Assume that $\psi\function{V(G)}{V(H)}$ is an
injection whose restriction to $V(G)\setminus D$ is the identity map.
Consider an arbitrary $k$-element set $U\subseteq V(G)$ and suppose
that $U=\{u_1,\ldots,u_a,v_1,\ldots,v_b\}$, where each
$u_j\in V(G)\setminus D$ and each $v_i\in D$. Then
$\psi(U)=\{u_1,\ldots,u_a,w_1,\ldots,w_b\}$, where $w_i=\psi(v_i)$.
Note that each $w_i\in D\cup(V(H)\setminus V(G))$. By A4 we have
$v_i\sim w_i$ in $H$ for all $i\le b$. By Item 3 of
Lemma \ref{lem:hsimprop2}, $\psi(U)\in E(H)$ iff $U\in E(H)$.
By A3, the latter is equivalent with $U\in E(G)$, proving that
$\psi$ is indeed a partial isomorphism from $G$ to $H$.

We now show that Conditions B1--B2 imply Conditions A1--A4.
For A1 and A2 this is trivial. Considering $\psi$ being the
identity map of $V(G)$ onto itself, we immediately obtain A3.
To obtain A4, it suffices to choose an arbitrary
$v_0\in[v]_G\cup(V(H)\setminus V(G))$ and prove that, for every
$v'\in [v]_G\cup(V(H)\setminus V(G))$, the vertices $v_0$ and $v'$
are similar in $H$. This will give
$[v]_G\cup(V(H)\setminus V(G))\subseteq[v]_H$.
The converse inclusion is given by Item 5 of Lemma \ref{lem:hsimprop2}.

Choose $v_0\in D$. According to
Item 1 of Lemma \ref{lem:hsimprop2}, we have to show,
for any $k$-vertex set $U\subseteq V(H)$, that $U\in E(H)$ iff
$U^{(v_0v')}\in E(H)$. If both $v_0$ and $v'$ or none of them are
in $U$, this is obvious. Otherwise it is enough to consider the case
that $v_0$ is in $U$ but $v'$ is not. Suppose that
$$
U=\{u_1,\ldots,u_a,v_0,v_1,\ldots,v_b,w_1,\ldots,w_c\},
$$
where all $u_i\in V(G)\setminus D$, all $v_i\in D$, and all
$w_i\in V(H)\setminus V(G)$.

Since $|D|\ge k$, we can choose in $D\setminus\{v_0,v_1,\ldots,v_b\}$
pairwise distinct vertices $v'_1,\ldots,v'_c$.
Define $\psi_1\function{V(G)}{V(H)}$ so that $\psi_1(v'_i)=w_i$
for each $i\le c$ and $\psi_1(x)=x$ for all other $x\in V(G)$.
If $v'\in[v]_G\setminus\{v'_1,\ldots,v'_c\}$, let $\psi_2$ be
the same as $\psi_1$. Notice that in this case
$\psi_2^{-1}(U^{(v_0v')})=(\psi_1^{-1}(U))^{(v_0v')}$
and $v_0\sim v'$ in $G$.
If $v'\in\{v'_1,\ldots,v'_c\}\cup(V(H)\setminus V(G))$,
let $\psi_2$ coincide with $\psi_1$ everywhere but $v_0$,
where we set $\psi_2(v_0)=v'$. Notice that now
$\psi_2^{-1}(U^{(v_0v')})=\psi_1^{-1}(U)$.
In both of the cases we have
$\psi_1^{-1}(U)\in E(G)$ iff $\psi_2^{-1}(U^{(v_0v')})\in E(G)$.
By B2, $\psi_1$ and $\psi_2$ are partial isomorphisms from $G$ to $H$.
It follows that $U\in E(H)$ iff $U^{(v_0v')}\in E(H)$,
completing derivation of~A4.

2)
Suppose that $H\cong K$ and $K\doteq G\add lv$.
Then Item 2 easily follows from Item 1 applied for $k$-graphs $G$ and $K$.
In particular, if we have Condition B2 with a set $D$, then
a map $\psi_0$ in Condition C2 can be taken the restriction of an
isomorphism from $K$ to $H$ to the set $V(G)\setminus D$.
\end{proof}

\begin{lemma}\label{lem:hd1}
Let $k\ge2$. If $G$ and $G'$ are non-isomorphic $k$-graphs of orders
$n$ and $n'$ respectively and $n\le n'$, then
\begin{equation}\label{eq:hd1n4}
\DD1{G,G'}\le\of{1-\frac1k}n+2k-1
\end{equation}
unless $G'=G\add(n'-n)v$ for some $v\in V(G)$.
\end{lemma}

\begin{proof}
We will describe a strategy of Spoiler winning $\game_r(G,G')$ for
$r=\lfloor(1-1/k)n+2k-1\rfloor$ unless $G'=G\add(n'-n)v$.
The strategy splits the game in two phases.

\smallskip

\centerline{\sc Phase 1}

\smallskip

Spoiler's aim is to select in $G$ a set of vertices $X$ with some
useful properties. He proceeds as follows. Initially, $X=\emptyset$.
If there is $B\subseteq\compl X$ with at most $k-1$ element such that
$\tau(X\cup B)>\tau(X)$, Spoiler chooses a such $B$ arbitrarily,
selects all vertices in $B$, and resets $X$ to $X\cup B$.
As soon as there is no such $B$, Phase 1 ends.

Suppose that Phase~1 has now ended, during which Spoiler made $s$
moves. The set $X=\{x_1,\ldots,x_s\}$ is from now on fixed and
consists of the vertices selected by Spoiler during Phase 1, where
$x_i$ is selected in the $i$-th round. It is easy to see that
\begin{equation}\label{eq:tau1}
\tau(X)-1\ge\frac s{k-1}.
\end{equation}
Since $s+\tau(X)\le n$, we have
\begin{equation}\label{eq:tau2}
s\le\of{1-\frac1k}(n-1).
\end{equation}
We will refer to the sets $Y(X)$ and $Z(X)$, and to the relation ${\simeq}_X$
omitting the subscript~$X$.

\begin{claim}\label{cl:simeq}
Let $U=\{u_1,\ldots,u_k\}$ and $W=\{w_1,\ldots,w_k\}$ be
$k$-element subsets of $V(G)$ and $u_i\simeq w_i$ for all $1\le i\le k$.
Then $U\in E(G)$ iff $W\in E(G)$.

In particular, if $u\simeq w$, then $u\sim w$ for any $u,w\in V(G)$,
i.e., $\calD(X)$ coincides with the partition of $Z$ into
the similarity classes.
\end{claim}

\begin{subproof}
Notice that the claim easily follows from its particular case that
$u_i=w_i$ for $i\le k-1$. We hence assume this.

Suppose on the contrary that $U\in E(G)$ but $W\notin E(G)$.
Let us modify $X$ in $k-1$ steps as follows. We will denote the result
of modification after the $i$-th step by $X_i$. Initially $X_0=X$.
We set $X_i=X_{i-1}\cup\{u_i\}$ if $u_i\notin X_{i-1}\cup Y(X_{i-1})$
and $X_i=X_{i-1}$ otherwise. We eventually enforce $u_k\notxyxk w_k$.
Since in the $i$-th step no single-element
${\equiv}_{X_{i-1}}$-class disappears, we have
$\tau(X_{k-1})>\tau(X)$ contradicting the assumption that Phase 1
has already ended.
\end{subproof}

Let $x'_i$ denote the vertex of $G'$ selected by Duplicator in the $i$-th
round of Phase 1 and $X'=\{x'_1,\ldots,x'_s\}$. Assume that Duplicator
has still not lost. Thus, the map $\phi\function X{X'}$ given by
$\phi(x_i)=x'_i$ for $i\le s$ is a partial isomorphism from $G$ to $G'$.
Similarly to Claim \ref{cl:phi2} we see that Duplicator now plays
under the following constraint.

\begin{claim}\label{cl:hphi2}
Whenever after Phase 1 Spoiler selects a vertex $v\in V(G)\cup V(G')$,
Duplicator responds with a $\phi$-similar vertex or otherwise
immediately loses.\qed
\end{claim}

\smallskip

\centerline{\sc Phase 2}

\smallskip

Similarly to the proof of Lemma \ref{lem:d1}, we conclude from
Claim \ref{cl:hphi2} that the $\phi$-similarity determines
the perfect matching between the singletons in $\calC(X)$
and the singletons in $\calC(X')$ unless Spoiler wins making 2 moves
and at most 1 alternation in Phase 2.
Denote the classes of $\calC(X)$ by $C_1,\ldots,C_t$
and the classes of $\calC(X')$ by $C'_1,\ldots,C'_{t'}$.
We therefore will assume that, for some $q\le t$, $|C_i|=1$ iff $i\le q$,
$|C'_i|=1$ iff $i\le q$, and $C_i\phieq C'_i$ for all $i\le q$.
Let $Y'=Y(X')$ and $\phi^*\function Y{Y'}$ be an extension of $\phi$
that maps $C_i$ onto $C'_i$ for every $i\le q$.

\begin{claim}\label{cl:hphiiso}
$\phi^*$ is a partial isomorphism from $G$ to $G'$, unless Spoiler
wins in the next $k$ moves with no alternation,
having made at total $s+k\le(1-1/k)n+k$ moves.
\end{claim}

\begin{subproof}
Assume that there is a $k$-vertex set $U\subseteq X\cup Y$
such that exactly one of the sets $U$ and $U'=\phi^*(U)$
is an edge of $G$ or $G'$ respectively. Let Spoiler select
all vertices of $U$. If Duplicator responds with the vertices
of $U'$, he obviously loses. If Duplicator moves at least once
outside $U'$, he violates the $\phi$-similarity and loses by
Claim~\ref{cl:hphi2}.
\end{subproof}

Assume that Duplicator is lucky to ensure that $\phi^*$ is a partial
isomorphism from $G$ to $G'$.
Let $Z=Z(X)$ and $Z'=Z(X')$.

\begin{claim}\label{cl:hphi*}
Whenever in Phase 2 Spoiler selects a vertex $v\in Z\cup Z'$,
Duplicator responds with a $\phi^*$-similar vertex or otherwise
loses in at most $k-1$ next rounds
with no alternation between $G$ and $G'$ in these rounds.
\end{claim}

\begin{subproof}
Let $u$ be the vertex selected by Duplicator in response to $v$
and assume that $u\notphistar v$. Suppose that $v\in Z'$
(the case of $v\in Z$ is completely similar). If $u\notin Z$,
Duplicator has already lost. If $u\in Z'$, there exists a
$(k-1)$-vertex set $W\subseteq X\cup Y$ such that $W\cup\{u\}$
is an edge but $\phi^*(W)\cup\{v\}$ is not or vice versa.
Spoiler selects the so far unselected vertices of $\phi^*(W)$.
Duplicator, who either responds with vertices in $W$ or
violates the $\phi$-similarity, loses.
\end{subproof}

Claim \ref{cl:hphi*} readily implies that every class in $\calD(X)$ or
$\calD(X')$ has a $\phi^*$-similar counterpart in, respectively,
$\calD(X')$ or $\calD(X)$ unless Spoiler wins making in Phase 2
at most $k$ moves and at most one alternation between the graphs
(selecting a vertex in a $\calD$-class without $\phi^*$-similar
counterpart and applying the strategy of Claim \ref{cl:hphi*}).
We will therefore assume that $\calD(X)=\{D_1,\ldots,D_p\}$,
$\calD(X')=\{D'_1,\ldots,D'_{p}\}$, and
$D_i\phistareq D'_i$ for all $i\le p$.

\begin{claim}\label{cl:hdidj}
Unless Spoiler is able to win making $2k-1$ moves and at most 1 alternation
in Phase 2, the following condition is met:
\begin{description}
\item[$(*)$]
For any $U=\{u_1,\ldots,u_k\}$, a $k$-vertex subset of $V(G)$,
and $U'=\{u'_1,\ldots,u'_k\}$, a $k$-vertex subset of $V(G')$,
such that $u_i\phistareq u'_i$ for all $i\le k$, we have
$U\in E(G)$ iff $U'\in E(G')$.
\end{description}
\end{claim}

\begin{subproof}
Suppose, to the contrary, that there are such $U$ and $U'$
but exactly one of $U$ and $U'$ is an edge.
Spoiler selects the vertices in $U'$. Let $U''=\{u''_1,\ldots,u''_k\}$
be the set of the respective responses of Duplicator.
If $u''_i\notphistar u'_i$ for some $i\le k$,
Spoiler wins in the next $k-1$ moves according
to Claim \ref{cl:hphi*}. Assume that $u''_i\phistareq u'_i$ for
all $i\le k$. Together with $u_i\phistareq u'_i$ for all $i$,
this implies $u_i\simeq u''_i$ for all $i$. Using Claim \ref{cl:simeq},
we conclude that $U''$ is an edge iff $U$ is and iff $U'$ is not.
Thus, Duplicator loses anyway.
\end{subproof}

In the sequel we suppose that the condition
$(*)$ in Claim \ref{cl:hdidj} holds (for else Spoiler wins
within the claimed bound \refeq{eq:hd1n4} for the number of moves).
Assume for a while that $|D_i|=|D'_i|$ for all $i\le p$.
Consider an arbitrary bijection $\bar\phi\function{V(G)}{V(G')}$
that extends $\phi^*$ and maps each $D_i$ onto $D'_i$.
The condition $(*)$ immediately implies that $\bar\phi$
is an isomorphism from $G$ to $G'$.
Since $G$ and $G'$ are supposed non-isomorphic, there must exist a $D_i$
such that $|D_i|\ne|D'_i|$. We will call such a $D_i$ {\em useful}.

Similarly to Claim \ref{cl:useful}, we obtain, as a corollary from
Claim \ref{cl:hphi*}, the following threat for Duplicator.

\begin{claim}\label{cl:huseful}
If $D_i$ is useful, then Spoiler is able to win
having made in Phase 2 at most $\min\{|D_i|,|D'_i|\}+k$ moves
and at most 1 alternation between $G$ and $G'$.\qed
\end{claim}

Suppose now that there are two useful classes, $D_i$ and $D_j$. Observe that
$$
|D_i|+|D_j|=|Z|-\sum_{l\ne i,j}|D_l|\le (n-s-q)-(\tau(X)-q-2)=
n-s-\tau(X)+2.
$$
It follows that one of the useful classes has at most $(n-s-\tau(X)+2)/2$
vertices. Thus, Spoiler is able to win totally in at most
\begin{equation}\label{eq:total}
s+\frac{n-s-\tau(X)+2}2+k=\frac{n+s-\tau(X)+2}2+k
\end{equation}
rounds. From \refeq{eq:tau1} and \refeq{eq:tau2}, we infer that
$$
s-\tau(X)\le s-\frac s{k-1}-1=s\,\frac{k-2}{k-1}-1\le\of{1-\frac 2k}(n-1)-1.
$$
The bound \refeq{eq:total} therefore does not exceed
$$
\of{1-\frac1k}n+k+\frac1k,
$$
which is within the required bound~\refeq{eq:hd1n4}.

Finally, suppose that there is a unique useful class $D_m$.
According to Claim \ref{cl:huseful}, Spoiler is able to win
in at most $|D_m|+k$ moves, with the total number of moves
$s+|D_m|+k$ within the required bound \refeq{eq:hd1n4}
provided $|D_m|\le k-1$.
Thus, we arrive at the conclusion that the bound \refeq{eq:hd1n4}
may not hold true in the only case that there is exactly one
useful class $D_m$ and $|D_m|\ge k$.
Let $\psi_0\function{V(G)\setminus D_m}{V(G')\setminus D'_m}$
be an extension of $\phi^*$ that maps each $D_i$, $i\ne m$, onto
$D'_i$. Let $\psi\function{V(G)}{V(G')}$ be an arbitrary injective
extension of $\psi_0$ (mapping $D_m$ into $D'_m$).
By the condition $(*)$ in Claim \ref{cl:hdidj}, $\psi$ is
a partial isomorphism from $G$ to $G'$.
Take an arbitrary $v\in D_m$. By Claim \ref{cl:simeq},
we have $D_m\subseteq[v]_G$. On the account of
Item 2 of Lemma \ref{lem:addequi} we conclude that
$G'=G\add(n'-n)v$, completing the proof of the lemma.
\end{proof}

\subsection{Defining $k$-graphs}

Theorem \ref{thm:def1} carries over $k$-graphs in a weaker form.
The analogs of Lemmas \ref{lem:sim} and \ref{lem:close} that we have
for $k$-graphs (see Remark \ref{rem:arbstr}) together with
Lemma \ref{lem:hd1}, give the following analog of
Lemma \ref{lem:upper}:

{\sl
Let $G$ and $G'$ be non-isomorphic $k$-graphs. Let $n$ denote
the order of $G$. Then we have
$$
\DD1{G,G'}\le\of{1-\frac1k}n+2k-1
$$
unless
\begin{equation}\label{eq:hs}
\gs(G)>\of{1-\frac1k}n+k-1.
\end{equation}
In the latter case we have
\begin{equation}\label{eq:hs12}
\gs(G)+1\le\DD1{G,G'}\le \gs(G)+k.
\end{equation}
}
Unfortunately, we cannot efficiently find the precise value of $\DD1{G,G'}$
in the range \refeq{eq:hs12} for $G$ with \refeq{eq:hs}, as it was done
for 2-graphs in Lemma \ref{lem:exact}.
Anyway, we have the following result, which is rather reasonable
as $\gs(G)$ is efficiently computable.

\begin{theorem}
If $\gs(G)\le\of{1-\frac1k}n+k-1$, then
$$
\DD1{G,G'}\le\of{1-\frac1k}n+2k-1.
$$
Otherwise,
$$
\DD1G\le\gs(G)+k
$$
and this bound is at most $k-1$ apart from the precise value of~$\D G$.
\end{theorem}

\section{Open questions}\label{s:open}
\mbox{}

\que\label{que:I}
In Section \ref{s:intro} we define
$$
\Ii n=\max\setdef{\D{G,G'}}{G\not\cong G',\ |G|=|G'|=n}.
$$
By Example \ref{ex:lower} and Theorem \ref{thm:d},
$$
\frac{n+1}2\le\Ii n\le\frac{n+3}2.
$$
This determines $\Ii n$ if $n$ is even and leaves two possibilites
for $\Ii n$ if $n$ is odd. Which is the right value?

\que
Given a graph $G$, is the number $\D G$ computable
(T.~\L uczak \cite{Luc})?
Can one, at least, improve the computable upper bound of
Theorem \ref{thm:def}, that is, can one lower the bound in
Theorem \ref{thm:def} below $n/2$, of course, extending the class~$\cl$?

\que
Theorem \ref{thm:D0} is an improvement in the alternation number
over Theorem \ref{thm:d}. Can one as well improve on Theorem \ref{thm:def}?

\que
Prove analogs of Theorem \ref{thm:D0} for digraphs and $k$-graphs.

\que
Improve the constant $c_d$ in Theorem \ref{thm:bdeg}.
Find connected bounded degree graphs $G$ and $G'$ of order $n$
with linear lower bound $\D{G,G'}=\Omega(n)$ better than that
given by Cai, F\"urer, and Immerman \cite{CFI}
(cf.\ Proposition \ref{prop:connbdeg}).

\que
For the optimum dimension of the \WL\/ algorithm we know bounds
$0.00465\,n<\Wl(n)\le 0.5\,n+0.5$, where the lower bound is due to
Cai, F\"urer, and Immerman \cite{CFI} and the upper bound is
shown in Section \ref{s:wl}. Make the gap between these bounds closer.

\que
Generalize Theorem \ref{thm:did1} over structures with maximum
relation arity $k$ proving a tight upper bound
$\DD1{G,G'}\le c_k n (1+o(1))$ for such structures
$G$ and $G'$ of order $n$. We currently know \cite{PVe} that
$1/2\le c_k\le 1-1/(2k)$. In the particular case of $k$-graphs,
Theorem \ref{thm:hd} gives us a better upper bound with $c_k=1-1/k$.
How tight is it?

\que
Find an analog of Lemma \ref{lem:exact} for $k$-graphs.
Namely, let $G$ and $G'$ be $k$-graphs of orders $n\le n'$,
$\gs(G)\ge n/2$, and $G'=G\add(n'-n)v$ for some $v\in V(G)$
with $\gs_G(v)=\gs(G)$. By an analog of Lemma \ref{lem:close}
for $k$-graphs (see Remark \ref{rem:arbstr}),
we have $\gs(G)+1\le\D{G,G'}\le\DD1{G,G'}\le \gs(G)+k$.
How to efficiently compute the exact value of $\D{G,G'}$?

\que
Let $G$ be a graph.
Let $L'(G)$ denote the minimum $l$ such that over
the variable set $\{x_1,\ldots,x_l\}$ there is a
first order formula defining $G$. Is it true that
$L'(G)=\max\setdef{\V{G,G'}}{G'\not\cong G}$?
In other terms, are the definabilities in the $l$-variable first order
logic and in the $l$-variable infinitary logic equivalent?

The affirmative answer would follow from this assumption:
There is a function $f$ such that, if graphs $G$ and $G'$ are
distinguished by a formula with $l$ variables, then they are
distinguished by a formula with $l$ variables of quantifier rank
at most $f(n,l)$, where $n$ is the order of $G$ (no dependence
on the order of $G'$!). In combinatorial terms:
If Spoiler can win the game on $G$ and $G'$ with $l$ pebbles in arbitrary
number of rounds (reusing pebbles is allowed), then he can win
with $l$ pebbles in  $f(n,l)$  rounds, irrespective of the order of $G'$.
Is it true?

\subsection*{Acknowledgments}
The third author is grateful to Georg Gottlob and members of his
department at the Vienna University of Technology for their hospitality
during this author's visit in the summer 2002.


\begin{thebibliography}{10}

\bibitem{AST}
N.~Alon, P.~Seymour, R.~Thomas.
\newblock
A separator theorem for nonplanar graphs.
\newblock
{\em J.\ Am.\ Math.\ Soc.} 3(4):801--808 (1990).

\bibitem{Bab}
L.~Babai.
\newblock
Automorphism groups, isomorphism, reconstruction.
\newblock
Chapter 27 of the {\em Handbook of Combinatorics,} 1447--1540.
R.~L.~Graham, M.~Gr\"otschel, L.~Lov\'asz Eds.
Elsevier Publ.\ (1995).

\bibitem{Bab2}
L.~Babai.
\newblock
On the complexity of canonical labelling of strongly regular graphs.
\newblock
{\em SIAM J.\ Comput.} 9:212--216 (1980).

\bibitem{BES}
L.~Babai, P.~Erd\H{o}s, S.~M.~Selkow.
\newblock
Random graph isomorphism.
\newblock
{\em SIAM J.\ Comput.} 9:628--635 (1980).

\bibitem{BKu}
L.~Babai, L.~Ku\v{c}era.
\newblock
Canonical labelling of graphs in linear average time.
\newblock
{\em Proc.\ of 20th Ann.\ Symp.\ on Foundations of Computer Science}
39--46 (1979).

\bibitem{Bol}
B.~Bollob\'as.
\newblock
The isoperimetric number of random regular graphs.
\newblock
{\em European J.\ Combin.} 9(3):241--244 (1988).

\bibitem{Bus}
P.~Buser.
\newblock
On the bipartition of graphs.
\newblock
{\em Discrete Appl.\ Math.} 9(1):105--109 (1984).

\bibitem{CFI}
J.-Y.~Cai, M.~F\"urer, N.~Immerman.
\newblock
An optimal lower bound on the
number of variables for graph identification.
\newblock
{\em Combinatorica} 12(4):389--410 (1992).

\bibitem{EFl}
H.-D.~Ebbinghaus, J.~Flum.
\newblock
{\em Finite model theory.}
\newblock
Springer Verlag, 2nd rev.\ ed.\ (1999).

\bibitem{Ehr}
A.~Ehrenfeucht.
\newblock
An application of games to the completeness problem for formalized
theories.
\newblock
{\em Fundam.\ Math.} 49:129--141 (1961).

\bibitem{Fra}
R.~Fra\"\i ss\'e.
\newblock
Sur quelques classifications des systems de relations.
\newblock
{\em Publ.\ Sci.\ Univ.\ Alger} 1:35--182 (1954).

\bibitem{Gro1}
M.~Grohe.
\newblock
Fixed-point logics on planar graphs.
\newblock
In: {\em Proc.\ of the Ann.\ Conf.\ on Logic in Computer Science}
6--15 (1998).

\bibitem{Gro2}
M.~Grohe.
\newblock
Isomorphism testing for embeddable graphs through definability.
\newblock
In: {\em Proc.\ of the 32nd ACM Ann.\ Symp.\ on Theory of
Computing (STOC)} 63--72 (2000).

\bibitem{GMa}
M.~Grohe, J.~Marino.
\newblock
Definability and descriptive complexity on databases of bounded tree-width.
\newblock
In: {\em Proc.\ of the 7th International Conference on
Database Theory}, Lecture Notes in Computer Science 1540, Springer-Verlag,
70--82 (1999).

\bibitem{Imm}
N.~Immerman.
\newblock
{\em Descriptive complexity.}
\newblock
Springer-Verlag (1999).

\bibitem{IKo}
N.~Immerman, D.~Kozen.
\newblock
Definability with bounded number of bound variables.
\newblock
{\em Information and Computation\/} 83:121--139 (1989).

\bibitem{ILa}
N.~Immerman, E.~Lander.
\newblock
Describing graphs: a first-order approach to graph canonization.
\newblock
In: {\em Complexity Theory Retrospective}, A.~Selman Ed.,
Springer-Verlag, 59--81 (1990).

\bibitem{KPSV}
J.~H.~Kim, O.~Pikhurko, J.~Spencer, O.~Verbitsky.
\newblock
How complex are random graphs in first order logic? In preparation (2003).

\bibitem{KMe1}
A.~V.~Kostochka, L.~S.~Melnikov.
\newblock
On bounds of the bisection width of cubic graphs.
\newblock
{\em Fourth Czechoslovakian Symposium on Combinatorics, Graphs
and Complexity} (Prachatice, 1990) 151--154.
\newblock
In: {\em Ann.\ Discrete Math.} 51, North-Holland, Amsterdam (1992).

\bibitem{KMe2}
A.~V.~Kostochka, L.~S.~Melnikov.
\newblock
On a lower bound for the isoperimetric number of cubic graphs.
\newblock
{\em Probabilistic methods in discrete mathematics}
(Petrozavodsk, 1992) 251--265.
\newblock
In: {\em Progr.\ Pure Appl.\ Discrete Math.} 1, Utrecht (1993).

\bibitem{Kuc}
L.~Ku\v{c}era.
\newblock
Canonical labelling of regular graphs in linear average time.
\newblock
{\em Proc.\ of 28th Ann.\ Symp.\ on Foundations of Computer Science}
271--279 (1987).

\bibitem{Luc}
T.~\L uczak.
\newblock
The 11-th Int.\ Conf.\ on Random Structures and Algorithm,
Poznan, 9--13 August 2003.

\bibitem{Pez}
E.~Pezzoli.
\newblock
Computational complexity of Ehrenfeucht-Fra\"\i ss\'e games on finite
structures.
\newblock
In: {\em Proc.\ of the CSL'98 Conf.},
G.~Gottlob, K.~Seyr Eds.
Lecture Notes in Computer Science 1584,
Springer-Verlag, 159--170 (1999).


\bibitem{PSV}
O.~Pikhurko, J.~Spencer, O.~Verbitsky.
\newblock
Succinct definitions in first order graph theory. In preparation (2003).

\bibitem{PVe}
O.~Pikhurko, O.~Verbitsky.
\newblock
Descriptive complexity of finite structures: saving the quantifier rank.
\newblock
Available at {\tt http://arxiv.org/abs/math.LO/0305244}

\bibitem{Spe}
J.~Spencer.
\newblock
{\em The strange logic of random graphs.}
\newblock
Springer Verlag (2001).

\bibitem{SJo}
J.~Spencer, K.~St.~John.
\newblock
The tenacity of Zero-One Laws.
\newblock
{\em The Electronic Journal of Combinatorics\/} 8(2) \#R17 (2001).

\bibitem{Tra}
B.~Trakhtenbrot.
\newblock
The impossibility of an algorithm for the decision problem for finite
models.
\newblock
{\em Dokl.\ Akad.\ Nauk SSSR} 70:596--572 (1950).
\newblock
English translation in: {\em AMS Transl.\ Ser.\ 2\/} 23:1--6 (1963).

\bibitem{Vau}
R.~Vaught.
\newblock
Sentences true in all constructive models.
\newblock
{\em Journal of Symbolic Logic\/} 25:39--58 (1960).

\bibitem{Ver}
O.~Verbitsky.
\newblock
The first order definability of graphs with separators
via the Ehrenfeucht game.
\newblock
Submitted (2003).

\end{thebibliography}
\end{document}